\documentclass[a4paper,12pt, psamsfonts, reqno]{amsart}
\usepackage{verbatim}

\usepackage{amsmath}
\usepackage{amssymb}
\usepackage{textcomp}
\usepackage{amsthm}
\usepackage{amsfonts}
\usepackage[all]{xy}
\usepackage{mathrsfs}
 \usepackage{graphicx}
 \usepackage{epstopdf}
 \usepackage{float}
 \usepackage{graphicx}
 
 \usepackage{epstopdf}
 \usepackage[T1]{fontenc}
 \usepackage{hyperref}
 
\newtheorem{theo}{Th\'eor\'eme}[section]

\def\Ind{{\mathop{\mathrm{Ind}}\nolimits}}

\def\C{\mathbb{C}}
\def\lC{\mathscr{C}}
\def\D{\mathscr{D}}

\def\N{\mathbb{N}}

\def\Q{\mathbb{Q}}

\def\Z{\mathbb{Z}}

\def\a{\mathbf{a}}
\def\e{\mathbf{e}}
\def\c{\mathbf{c}}
\def\m{\mathfrak{m}}
\def\b{\mathbf{b}}
\def\d{\mathbf{d}}

\def\O{\mathcal{O}}
\def\line{\overline}

\def\Id{\mathop{\mathrm{Id}}\nolimits}
\def\Hom{\mathop{\mathrm{Hom}}\nolimits}

\def\ker{\mathop{\mathrm{ker}}\nolimits}

\def\Im{\mathop{\mathrm{Im}}\nolimits}

\def\dim{\mathop{\mathrm{dim}}\nolimits}
\def\rank{\mathop{\mathrm{rank}}\nolimits}

\def\sym{\mathop{\mathrm{sym}}\nolimits}

\def\remk{\noindent\textit{Remark:~}}

\newtheorem{prop}[theo]{Proposition}
\newtheorem{prop-def}[theo]{Proposition-Definition}
\newtheorem{def-prop}[theo]{Definition-Proposition}
\newtheorem{cor}[theo]{Corollary}

\newtheorem{lemma}[theo]{Lemma}
\newtheorem{teo}[theo]{Theorem}
\newtheorem{definition}[theo]{Definition}

\newtheorem{conj}[theo]{Conjecture}

\newtheorem{notation}[theo]{Notation}
\newtheorem{example}[theo]{Example}

\allowdisplaybreaks

\everymath{\displaystyle}

\input xy
\xyoption{all}

\title{Symmetrization of representations of $GL_N$}

\author{Taiwang, DENG}


\address{Max Planc Institut f{\"{u}}r Mathematik, 
 Vivatsgasse 7, 53111 Bonn.}
\date{}

\begin{document}

\begin{abstract}
In this article, we develop a process to
symmetrize the irreducible admissible representation
of $GL_N(\Q_p)$, as a consequence we obtain a more geometric
understanding of the coefficient $m(\b, \a)$ appearing 
in the decomposition of parabolic inductions, which allows
us to prove a conjecture posed by Zelevinsky.
\end{abstract}

\maketitle

\tableofcontents

\section{Introduction}

In this article, we will develop a process to symmetrize the admissible 
representations of general linear groups over non-archimedean fields.

For a $p$-adic field $F$ and $g>1$, an irreducible admissible representation $\rho$
of $GL_{g}(F)$ is called cuspidal if for all proper parabolic
subgroup $P$, the corresponding Jacquet functor $J^G_{P}$ sends $\rho$ to $0$. We 
write 
$$\nu: GL_g(F)\rightarrow \C, \qquad\nu(x)=|\det(x)|$$ 
and for $k\geq 1$ and $\rho$ a cuspidal irreducible representation of $GL_g(F)$, 
we call the set $$\Delta_{\rho, k}=\{\rho, \rho\nu, \cdots, \rho\nu^{k-1}\}$$
a segment. For such a segment, the normalized induction functor
$$\Ind_{P_{g,\cdots,g}}^{GL_{kg}(F)} (\rho \otimes \cdots \otimes \rho\nu^{k-1})$$
contains a unique irreducible sub-representation denoted by $L_{[\rho,\nu^{k-1}\rho]}$, where
$P_{g,\cdots,g}$ is the standard parabolic subgroup with Levi subgroup isomorphic to $k$ blocks of $GL_g$.
Then a multisegment is a multiset of segments, by multiset we mean a set with multiplicities.
For $i=1, \cdots, r$, let $\rho_i$ be an irreducible cuspidal representation of $GL_{n_i}(F)$ and for $k_i\in \N$, 
by definition, the multisegment
\[
 \a=\{\Delta_{\rho_i, k_i}: i=1, \cdots, r\}, 
\]
is of degree $\deg(\a)=\sum n_ik_i$.
In \cite{Z2}, the author gave a parametrization $\a\mapsto L_{\a}$
of irreducible admissible representations of $GL_{n}(F)$
in terms of multisegments of degree $n$, where for a multisegment $\a$
with suitable order on its elements(cf. Theorem \ref{teo:1}),
the representation
$L_{\a}$
is the unique irreducible submodule of the parabolic induced representation
\[
 \pi(\a)=\Ind_{P}^{GL_{n}(F)}(L_{\Delta_{\rho_1, k_1}}\otimes \cdots \otimes L_{\Delta_{\rho_r, k_r}}).
\]
Now given two multisegments $\a$
and $\b$, one wants to determine the multiplicity 
$m(\b, \a)$ of $L_{\b}$ in $\pi(\a)$.

Thanks to the Bernstein central decomposition, one is  reduced to the case 
where the cuspidal supports of $\a$ and $\b$ belong to the 
same Zelevinsky line $\{\rho_0\nu^k: k \in \Z\}$. Zelevinsky also conjectured that 
$m(\b, \a)$ is independent of $\rho_0$ and 
depends only on the relative position of $\a$ and $\b$: this conjecture now follows from the theory of types,  
cf. \cite{MV}. So one is reduced to the simplest case where $\rho_0$ is the trivial representation.

\bigskip

Let us now explain what is known about these coefficients $m(\b,\a)$ where the cuspidal support of
$\a,\b$ belongs the Zelevinsky line of the trivial representation.
First of all, it is proved in \cite{Z2} that there exists a poset structure on 
the set of multisegments such that 
$m_{\b, \a}>0$ if and only if $\b\leq \a$. And we let 
\[
  S(\a)=\{\b: \b\leq \a\}.
 \]
In \cite{Z3},  Zelevinsky introduced the nilpotent
orbit associated to a multisegment $\a$. 
More precisely, to a multisegment $\a$, one can 
associate $\varphi_{\a}: \Z\rightarrow \N$ 
with $\varphi_{\a}(k)$ the multiplicities
of $\nu^k$ appearing in $\a$. For each $\varphi$, 
$V_{\varphi}$ is a $\C$-vector space of dimension 
$\deg \varphi: =\sum_{k\in \Z}\varphi(k)$ with 
graded $k$-part of dimension $\varphi(k)$. 
Then $E_{\varphi}$ is the set of endomorphisms 
$T$ of degree $+1$, which admits 
a natural action of the group $G_{\varphi}=\prod_{k}GL(V_{\varphi, k})$.
Then the orbits of $E_{\varphi}$ under $G_{\varphi}$
are parametrized by multisegments $\a=\sum_{i\leq j}a_{ij}\Delta_{\nu^i, j-i+1}$ such that 
$\varphi=\varphi_{\a}$ consists
of $T$ with $a_{ij}$ Jordan cells starting from $V_{\varphi, i}$ and ending in $V_{\varphi,j}$. 
We denote by $O_{\a}$ this orbit and we have the following nice property
\[
 \line{O}_{\a}=\bigsqcup_{\b\geq \a}O_{\b}.
\]
Now given a local system $\mathcal{L}_{\a}$
on $O_{\a}$, we can consider its 
intermediate extension $IC(\mathcal{L}_{\a})$ 
on $\line{O}_{\a}$ and its fiber at 
a geometric point $z_{\b}$ of $O_{\b}$. We form the 
Kazhdan-Lusztig polynomial 
\[
 P_{\a, \b}(q)=\sum_{i}q^{i/2}\dim_{\C}\mathcal{H}^i(IC(\mathcal{L}_{\a}))_{z_{\b}}.
\]
Zelevinsky then conjectured that $m_(\b, \a)=P_{\a, \b}(1)$ and 
call it the $p$-adic analogue of Kazhdan Lusztig Conjecture. 
This conjecture is a special case of a more general 
multiplicities formula proved by Chriss and Ginzburg in \cite{CG}, 
chapter 8.

\bigskip

\emph{In this paper}, we first introduce the notion of a symmetric multisegment (cf. definition \ref{def: 2.1.5}), which is, 
roughly speaking,  a multisegment such that the 
beginnings and the ends of its segments are distinct and its segments admit non-empty intersections. 
We show that for a well chosen\footnote{Thanks to corollary \ref{coro-sym1} which is a particular case
of the Zelevinsky's conjecture, the results are independent of the choice of $\a_{\Id}$.}
symmetric multisegment $\a_{\Id}$, there is a natural bijection 
between the symmetric group $S_n$ to the set of symmetric multisegments $S(\a_{\Id})$, cf.  proposition \ref{teo: 2.3.2},
where 
$n$ is the number of  segments contained in  $\a_{\Id}$.

When we restrict the geometry of the nilpotent orbits 
to the symmetric locus, 
%
%
we recover the  geometric situation of the Schubert varieties associated 
to $S_n$ and obtain that for two symmetric multisegment 
$\a_{\sigma}, \a_{\tau}$ associated to $\sigma, \tau\in S_n$, 
the coefficient $m(\a_{\sigma}, \a_{\tau})=P_{\sigma, \tau}(1)$. 

\bigskip

The next step \emph{in section 4} is to try to reach non symmetric 
cases, starting with symmetric ones. 
For example for $\a\geq \b$ two multisegments and $\nu^k$ in 
the cuspidal support of $\a$,  
one can eliminate every $\nu^k$ which appears at the end of 
some segments in $\a$ and $\b$ to obtain 
respectively a new pair of multisegments $\a^{(k)}, ~\b^{(k)}$ and 
try to prove that that $m(\b, \a)=m(\b^{(k)}, \a^{(k)})$. This 
result is almost true if we demand that $\b$ belongs to 
some subset $S(\a)_k$ of $S(\a)$, cf. Proposition\ref{cor: 3.2.3}.
The proof relies on the study of the geometry of nilpotent orbits and their links
with the Grassmannian, cf. the introduction of section 3.

%

\bigskip

In\emph{ section 5}, 
we prove the main result in the present paper, which is an iteration
of the process we develop in section 4.
In fact, for  a multisegment $\a$ and $k_1, \cdots, k_r$
integers such that $\nu^{k_i}$ appears in the 
supercuspidal support of $\a$, let 
\[
 \a^{(k_1, \cdots, k_r)}=(((\a^{(k_{1})})\cdots)^{ (k_{r})}), 
\]
and 
\[
S(\a)_{k_{1}, \cdots, k_{r}}=
\{\c\in S(\a):
\c^{( k_{1}, \cdots , k_{i})}
\in  S(\a^{( k_{1}, \cdots , k_{i})})_{k_{i+1}}, \text{ for }i=1, \cdots, r\}. 
\]
Then we show that for $\b\in S(\a)_{k_{1}, \cdots, k_{r}}$,
we always have 
\[
 m(\b, \a)=m(\a^{(k_1, \cdots, k_r)}, \b^{(k_1, \cdots, k_r)}),
\]
Reciprocally, we show,  cf. proposition \ref{prop: 4.2.4}, that for any pair of multisegments
$\a>\b$, we can find $\a^{\sym}$ and 
$\b^{\sym}<\a^{\sym}$ such that 
\[
m(\b, \a)=m(\b^{\sym}, \a^{\sym}).
\]
This finishes our symmetrization of a multisegment.As an application, we give a proof for a conjecture due to Zelevinsky.  

\begin{conj}(cf. \cite{Z2} \S 8.7)
The coefficient $m(\b, \a)$ when $\a$ and $\b$ belong
to the same $\O(\Pi)$, depends only on the mutual relation between $\a$ and $\b$.
\end{conj}

As a final remark, we mention the recent paper \cite{LM17}, where the authors
state an open orbit conjecture(cf. Conjecture 1.1), which is a variant
of the open orbit conjecture 
by Geiss-Leclerc-Schr\"{o}er \cite{GLS}.
They proved their conjecture in regular 
case through proving in p-adic setting the results 
of \cite{KKMO}, and left open the irregular cases.
Our symmetrization method seems to be a perfect 
way to deduce results for irregular multisegments 
from the regular ones. We pursue in this direction in 
a subsequent paper. To complete our introduction, 
we should also mention that a large part of the original open orbit conjecture
(those basis corresponding to cluster monomials in the sense of cluster algebras) is also known now by Kashiwara 
and his collaborators. 

\par \vskip 1pc
{\bf Acknowledgements}
This paper is part of my thesis at University Paris 13, which is funded by the program DIM of the
region Ile de France. I would like to thank my advisor Pascal Boyer for his keen interest in this work and 
his continuing support and countless advice.
It is rewritten during my first year of posdoc at Universit\"{a}t Bonn, I thank their hospitality. In addition, 
I would like to thank Alberto M\'inguez et Vincent S\'echerre, Yichao TIAN for their helpful discussions on the subject. 
Finally, I thank Bernard Leclerc and Erez Lapid for their careful reading and suggestions of my thesis.
\vspace{1cm}

\section{Zelevinsky classfication of induced representations}
In this section we recall Zelevinsky
classification of induced representations of $GL_N(F)$, with
$F$ a finite extension of $\mathbb{Q}_p$.

\begin{notation}
We fix a uniformizer $\varpi_F$ of $F$, and an absolute vaule 
$|.|$ on $F$ such that $|\varpi_F|=1/q$, where $q$ is the 
cardinal of its residue field. Note $\nu$ the character of $GL_n(F)$ defined by $\nu(g)=|\det{g}|$.
\end{notation}

\begin{definition}
By segment $\Delta$, we mean a finite consecutive subset of integers 
\[
\Delta=\{k_1, k_1+1, \cdots, k_2\}, \quad k_1\leq k_2, \quad (k_1, k_2)\in \Z^2.
\]
And we define a multisegment $\m$ to be a multiset of segments, 
\[
\m=\{\Delta_1, \cdots, \Delta_r\}.
\]
And we call 
\[
\deg{\a}=\sum_{i=1}^{r}\sharp{\Delta_i}
\]
the degree of $\a$.
\end{definition}

Following Zelevinsky, 

\begin{prop-def}
For any irreducible cuspidal representation $\rho$ of $GL_n(F)$ and a segment $\Delta$, we can associate 
an induced irreducible representation $L_{(\Delta,\rho)}$ of $GL_{n\deg{\Delta}}(F)$ in a 
unique way. When $\rho=1$ be the trivial character 
of $GL_1(F)$, we write directly $L_{\Delta}$.
\end{prop-def}

\begin{definition}
\begin{description}
\item[(1)]We say two segments $\Delta_1$ and $\Delta_2$ are 
linked if $\Delta_1\cup \Delta_2$ is again a segment and different from $\Delta_1$ and $\Delta_2$.
\item[(2)]We define the following partial order on 
the set of segments
\begin{displaymath} \left\{ \begin{array}{cc}
&[j, k]\prec  [m,n], \text{ if } k< n,\\
&[j,k]\prec  [m,n], \text{ if } j>m, n=k.
\end{array}\right. 
\end{displaymath}

\end{description}

\end{definition}

\begin{proof}
For explicite constuction, we refer to \cite{Z2}.
\end{proof}

\begin{definition}
For any pair of representation $(\pi_1, \pi_2)\in Rep(GL_n(F))\times Rep(GL_m(F))$, let $\pi_1\times \pi_2$ be 
the normalized induction of $\pi_1\otimes \pi_2$, which is a representation of
$GL_{n+m)}(F)$.

\end{definition}

\begin{prop}(\cite{Z2} Theorem 4.2)\label{prop: 1.1.11}
The following are equivalent:
\begin{description}
\item[(1)]The induced representation
\[
L_{(\Delta_1, \rho)}\times L_{(\Delta_2, \rho)}\times \cdots L_{(\Delta_r, \rho)}
\]
is irreducible.
\item[(2)]For any $1\leq i, j\leq r$, the segments
$\Delta_i$ and $\Delta_j$ are not linked with each other.
\end{description}
\end{prop}

\begin{prop}(cf. \cite{Z3} section 4.6 )\label{prop: 1.1.15}
 Let $\Delta_{1}$ and $\Delta_{2}$ be two linked segments with $\Delta_1\preceq \Delta_2$, then 
 \[
  L_{(\Delta_{1}, \rho)}\times L_{(\Delta_{2}, \rho)}
 \]
contains a unique sub-representation $L_{\a_1}$ and a unique quotient $L_{\a_2}$ with 
$$\a_{1}=\{(\Delta_{1}, \rho), ~(\Delta_{2}, \rho)\}, \quad \a_{2}=\{(\Delta_{1}\cup \Delta_{2}, \rho), ~(\Delta_{1}\cap \Delta_{2},\rho)\}.$$
\end{prop}
 
\begin{definition}
Let $\a=\{\Delta_1, \Delta_2, \cdots, \Delta_r\}$ be a multisegment
such that $\Delta_1$ and $\Delta_2$ are linked. By an elementary 
operation we mean replacing the segments $\Delta_1$ and $\Delta_2$
by $\Delta_1\cap \Delta_2$ and $\Delta_1\cup \Delta_2$. In this case,
we say $\a'=\{\Delta_1\cap \Delta_2, \Delta_1\cup \Delta_2, \cdots, \Delta_r\}$ is obtained from $\a$ via an elementary operation.
\end{definition} 
 
\begin{definition}
We define $\b\leq \a$  if $ \b$ can be obtained from $\a$ via 
a sequence of elementary operations. Denote
\[
S(\a)=\{\b: \b\leq \a\},
\]
then $\leq$ defines a partial order on $S(\a)$(cf. \cite{Z2} 7.1).
\end{definition}

We recall the following classifying theorem due to Zelevinsky.

\begin{teo}\label{teo:1}(\cite{Z2} Theorem 6.1)
Let $\a=\{(\Delta_1, \rho_1), \cdots, (\Delta_r, \rho_r)\}$ be 
a multisegment of cuspidal representations
with $\Delta_1\succeq \Delta_2\succeq \cdots \succeq \Delta_r$
, then 
\begin{description}
\item[(1)] The representation 
\[
L_{(\Delta_1, \rho_1)}\times \cdots \times L_{(\Delta_r, \rho_r)}
\]
contains a unique sub representation, which we denote by $L_{\a}$.

\item[(2)] The representation $L_{\a'}$ and $L_{\a}$ are isomorphic if and only if $\a=\a'$.
\item[(3)]Any irreducible representation of $GL_n(F)$ is of the 
form $L_{\a}$.

\end{description}

\end{teo}

\begin{definition}
A multisegment of cuspidal representations $$\a=\{(\Delta_1, \rho_1), \cdots, (\Delta_r, \rho_r)\}$$ 
is said to be well ordered if $\Delta_1\succeq \Delta_2\succeq \cdots \succeq \Delta_r$. 
\end{definition}

\begin{notation}
We denote by $\mathcal{R}_{n}$ the Grothendieck group of the category
of finite length representations of $GL_{n}(F)$ and
$$\mathcal{R}^{univ}=\oplus_{n\geq 1} \mathcal{R}_{n}.$$
\end{notation}

As was observed by Zelevinsky, the group 
$\mathcal{R}^{univ}$ can 
be endowed with a Hopf algebra structure via 

\begin{prop}
The set $\mathcal{R}^{univ}$ is a bi-algebra with the multiplication $\mu$ and co-multiplication $c$ given by
  \[
   \mu(\pi_{1}\otimes \pi_{2})=\pi_{1}\times \pi_{2},\qquad 
   c(\pi)=\sum_{r=0}^{n} J^{GL_{n}(F)}_{P_{r,n-r}}(\pi),
  \]
where $J^{GL_{n}(F)}_{P_{r,n-r}}$ denotes
the Jacquet functor from the category of smooth representations 
of $GL_n(F)$ to the category of smooth representations of $M_{r,n-r}=GL_r(F)\times GL_{n-r}$ regarded as the 
Levi subgroups of $P_{r, n-r}$, where $P_{r, n-r}$ is the 
unique parabolic subgroup containing the upper triangular matrices
with the given Levi subgroups.
\end{prop}

Now Zelevinsky's classification theorem can be reformulate into the following 

 \begin{cor}\label{cor: 1.2.3}
 The algebra $\mathcal{R}^{univ}$ is a polynomial ring with indeterminates $\{L_{\Delta}: \Delta\in \Sigma^{univ}\}$.
 Moreover, as a $\Z$-module,
  the set $\{L_{\a}: \a\in \O^{univ}\}$ form a basis for $\mathcal{R}^{univ}$.
 \end{cor}

 \remk
 Note that this implies the Bernstein Center theorem, i.e, we have a decomposition
 \[
  \mathcal{R}^{univ}=\prod_{\rho}\mathcal{R}(\rho),
 \]
where $\rho$ runs through the equivalent classes of irreducible (super)cuspidal representations.
Here we say two irreducible (super)cuspidal
representations $\rho, \rho'$ are equivalent if 
\[
\rho'\in \Pi_{\rho}=\{\rho \nu^s: s\in \Z\}.
\]
We denote by $\O(\rho)$ the set of multisegments supported on $\Pi_{\rho}$.

\begin{notation}\label{nota: 1.1.15}
From now on, for $\a=\{\Delta_{1}, \cdots, \Delta_{r}\}$ being well ordered,  we denote
\[
 \pi(\a)=L_{\Delta_{1}}\times \cdots \times L_{\Delta_{r}}.
\]

\end{notation}

 According to Theorem \ref{teo:1}, let $\a=\{\Delta_{1}, \cdots, \Delta_{r}\} $ be a multisegment
 with support contained in some Zelevinsky line $\Pi_{\rho}$, then  we can write
 \addtocounter{theo}{1}
 \begin{equation}\label{eq: m(b, a)}
  \pi(\a)=\sum_{\b\in \O(\rho)} m(\b,\a)L_{\b}
 \end{equation}
where $\pi(\a)=\Delta_{1}\times \cdots \times \Delta_{r},~ m(\b,\a) \in \mathbb N$. The aim of this paper is to give some new insights on these
$m(\b,\a)$.

\remk
it is conjectured in \cite{Z2} 8.7 that the coefficient $m(\b,\a)$ depends only on the combinatorial relations of 
$\b$ and $\a$, and not on the specific cuspidal representation $\rho$. The 
independence of specific cuspidal representation can be shown by type theory,
see for example \cite{MV}. \textbf{In 
other words, as far as we are concerned with the coefficient $m(\b,\a)$, we can restrict
ourselves to the special case $\rho=1$, the trivial representation of $GL_{1}(F)$.}

In final part of this section we show how to define some analogue of the Zelevinsky derivation, 
which serve as a tool for us in the sequel and motivates the development of this paper. 

\begin{definition}
 We define a left partial derivation with respect to index $i$ to be a morphism of algebras
 \begin{align*}
 &^i\D: \mathcal{R}\rightarrow \mathcal{R},\\
 &^i\D(L_{[j,k]})=L_{[j,k]}+\delta_{i,j}L_{[j+1,k]} \text{ if }(k>j),\\
 &^i\D(L_{[j]})=L_{[j]}+\delta_{[i],[j]}.
 \end{align*}
 Also we define a right partial derivation with respect to index $i$ to be a morphism of algebras
 \begin{align*}
 &\D^i: \mathcal{R}\rightarrow \mathcal{R}\\
 &\D^i(L_{[j,k]})=L_{[j,k]}+\delta_{i,k}L_{[j,k-1]} \text{ if }(j<k)\\
 &\D^i(L_{[j]})=L_{[j]}+\delta_{[i],[j]}.
 \end{align*}
\end{definition}

\begin{definition}
 We define 
 \[
  \D^{[i,j]}=\D^{j}\circ \cdots \circ \D^{i}
  \]
  \[
  ^{[i,j]}\D=(^{i}\D)\circ \cdots \circ ( ^{j}\D)
 \]
 And for $\c=\{\Delta_{1}, \cdots, \Delta_{s}\}$ with
 \[
  \Delta_{1}\preceq \cdots \preceq \Delta_{s},
 \]
we define 
\[
 \D^{\c}=\D^{\Delta_{1}}\circ \cdots \circ \D^{\Delta_{s}}
\]
and
\[
  ^{\c}\D=(^{\Delta_{s}}\D)\circ \cdots \circ (^{\Delta_{1}}\D).
\]
\end{definition}

\remk we recall that in \cite{Z1} 4.5, Zelevinsky defines a derivative $\mathscr{D}$ to be 
 an algebraic morphism 
 \[
  \mathscr{D}: \mathcal{R}\rightarrow \mathcal{R},
 \]
which plays a crucial role in Zelevinsky's classification theorem.

The relation between Jacquet functor and derivative is given by
\begin{prop}(cf. \cite{Z2}3.8)
 Let $\delta$ be the algebraic morphism such that $\delta(\rho)=1$ for all
 $\rho\in \lC$ and $\delta(L_\Delta)=0$ for all non cuspidal representations $L_\Delta$. Then
 \[
  \mathscr{D}=(1\otimes \delta)\circ c,
 \]
where $c$ is the co-multiplication.
\end{prop}

 The main advantage to work with partial derivatives instead of the derivative defined by
 Zelevinsky is that they are much more simpler but share the following positivity properties:
 
 \begin{teo}\label{teo: 3}
  Let $\a$ be any multisegment, then we have 
  $$\D^i(L_{\a})=\sum_{\b\in \O}n(\b,\a)L_{\b},$$
  such that $n(\b,\a)\geq 0$, for all $\b$.
 \end{teo}
 
\remk the same property of positivity holds for $^i\D$.

The theorem follows from the following two lemmas

\begin{definition}
For $i\in \Z$, let $\phi_{i}$ be the morphism of algebras defined by
\begin{align*}
 \phi_{i}: \mathcal{R}&\rightarrow \Z\\
 \phi_{i}([j,k])&=\delta_{[i], [j,k]}.
\end{align*}
\end{definition}

\begin{lemma}
For all multisegment $\a$, we have $\phi_{i}(L_{\a})=1$ if and
only if $\a$ contains no other segments than $[i]$, otherwise it is zero.
\end{lemma}

\begin{proof}
We prove this result by induction on the cardinality  of $S(\a)$, denoted by $|S(\a)|$. If $|S(\a)|=1$, then
$\a=\a_{min}$, hence $\phi_{i}(L_{\a})= \phi_{i}(\pi(\a))$, which is nonzero if and
only if $\a$ contains no other segments than $[i]$, and in latter case it is 1. 
Let $\a$ be a general multi-segment,
\[
 \pi(\a)= L_{\a}+\sum_{\b<\a}m(\b,\a)L_{\b}.
\]
Now $|S(\a)|>1$, we know that $\a$ is not minimal in $S(\a)$, hence $\a$
contains segments other than $[i]$, which implies $\phi_{i}(\pi(\a))=0$.

Since $|S(\b)|<|S(\a)|$ for any $\b<\a$, by induction, we know that $\phi_{i}(L_{\b})=0$
because $\b$ must contain segments other than $[i]$. So we are done.

\end{proof}

\begin{lemma}
 We have $\D^i= (1\otimes \phi_{i})\circ c$.
\end{lemma}

\begin{proof}
Since both are algebraic morphisms, we only need to check that they coincide on
generators. We recall the equation from \cite{Z2}, proposition 3.4
\[
 c(L_{[j,k]})=1\otimes L_{[j,k]}+\sum_{r=j}^{k-1}L_{[j,r]}\otimes L_{[r+1, k]}+L_{[j,k]}\otimes 1.
\]
   Now applying $\phi_{i}$,
\begin{align*}
 (1\otimes \phi_{i})c(L_{[j,k]})=&L_{[j,k]}+\delta_{i,k}L_{[j,k-1]} \text{ if }(k>j)\\
 (1\otimes \phi_{i})c(L_{[j]})=&L_{[j,k]}+\delta_{i,j},
 \end{align*}
 where $\delta_{i, j}$ is the Kronecker symbol. 
Comparing this with the definition of $\D^i$ yields the result.
\end{proof}

\remk We have the following relation between partial derivative and derivative of Zelevinsky.
 Let $e(\a)=\{[i_{1}], \cdots, [i_{\alpha}]: i_{1}\leq \cdots \leq i_{\alpha}\}$ be the end
of $\a$, then 
 \[
  \D(\a)=\D^{[i_{1}, i_{\alpha}]}(\a).
 \]

\section{Orbital varieties and KL polynomials}

A geometric interpretation of Zelevinsky's classification, which is also due to Zelevinsky, is 
to consider the orbital varieties associated to
multisegments, cf. \cite{Z3}, \cite{Z4}. 

\begin{definition}
Let $\a$ be a multisegment, we define a function
\[
\varphi_{\a}: \Z\rightarrow \N
\]
by letting $\a=\sum_{i\leq j}a_{ij}[i, j]$, and
\[
\varphi_{\a}(k)=\sum_{i\leq k\leq j}a_{ij}.
\]
We call $\varphi_{\a}$ the weight function of $\a$.
\end{definition}

\begin{definition}
Let $\varphi:\Z\rightarrow \N$ be a function with
finite support. Consider the $\Z$-graded $\C$-vector space
 \[
 V_{\varphi}=\oplus_{k}V_{\varphi, k}, \quad \dim(V_{\varphi, k})=\varphi(k).
 \]
Moreover,
 \begin{description}
  \item[(1)] let $E_{\varphi}$ be the set of endomorphisms of $V_{\varphi}$ of degree 1;
  \item[(2)] let $G_{\varphi}=\prod_{k}GL(V_{\varphi, k})$ be the automprhism of 
  $V_{\varphi}$.
 \end{description}
 
\end{definition}

\remk
The group $G_{\varphi}$ acts naturally on the endomorphism space $E_{\varphi}$
via conjugations.

\begin{prop}(cf.\cite{Z3}, 2.3)\label{prop: 2.2.4}
The orbits of $E_{\varphi}$ under
$G_{\varphi}$ are naturally parametrized by multisegments of weight $\varphi$. Moreover, let $O_{\a}$ be the orbit associated
to the multisegment $\a$, then 
\[
\a\leq \b\Leftrightarrow O_{\b}\subseteq \line{O}_{\a}.
\]
\end{prop}

\begin{proof}
Let $\a=\sum_{i\leq j}a_{ij}[i,j]$ such that 
$\varphi_{\a}=\varphi$, then 
the orbit associated consists of 
the operators having exactly $a_{ij}$ Jordan cells starting from $V_{\varphi, j}$ and
ending in $V_{\varphi,i}$.
\end{proof}

\begin{example}
We consider the function $\varphi:\Z\rightarrow \N$ with
\[
\varphi(0)=\varphi(1)=2, \quad \varphi(i)=0, \quad \forall i\neq 1, 2.
\]
Then $E_{\varphi}=\{T: V_{\varphi,0}\rightarrow V_{\varphi, 1}\}$.
In this case $E_{\varphi}$ contains 3 orbits which are determined by
the rank of $T$:
\begin{description}
\item[(1)]the orbit $\{T: \rank{T}=0\}=O_{\a_0}$ with $\a_0=\{[0], [0], [1], [1]\}$;
\item[(2)]the orbit $\{T: \rank{T}=1\}=O_{\a_1}$ with $\a_1=\{[0], [1], [0, 1]\}$;
\item[(3)]the orbit $\{T: \rank{T}=2\}=O_{\a_2}$ with $\a_2=\{[0, 1], [0,1]\}$.
\end{description}

\end{example}

\remk
The orbits $\{O_{\b}: \varphi_{\b}=\varphi\}$ give rise to a stratification of the affine space $E_{\varphi}$.

\begin{definition}
Let $\a$, $\b$ be two multisegments such that $\b\in S(\a)$. Then we define 
the polynomial
\[
 P_{\a, \b}(q)=\sum_{i}q^{(i+d_{\b})/2}\dim \mathcal{H}^{i}(\line{O}_{\b})_{x_{\a}},
\]
where 
\begin{itemize}
\item $\mathcal{H}^{i}(\line{O}_{\b}):=\mathcal{H}^{i}(IC(\line{O}_{\b}))$ is the intersection complex 
supported on $\line{O}_{\b}$ which is constant with stalk
$\C$ on $O_{\b}$;
\item $x_{\a}\in O_{\a}$ is an arbitrary point and $d_{\b}=\dim(O_{\b})$. 
\end{itemize}
We call it the Kazhdan Lusztig polynomial associated to $\{\a, \b\}$.
\end{definition}

We recall the following fundamental result, 
which is conjectured by Zelevinsky and named by whom the 
p-adic analogue of Kazhdan Lusztig conjecture.

\begin{teo}(\cite{Z3}, \cite{CG})\label{teo: 4.1.5}
 Let $\mathcal{H}^{i}(\line{O}_{\b})_{\a}$ denote the stalk at a point
 $x\in O_{\a}$ of the $i$-th intersection cohomology sheaf of the variety $\line{O}_{\b}$.
 Then 
 \[
  m(\b,\a)=P_{\b, \a}(1).
 \]
 \end{teo}

\remk In \cite{Z4} Theorem 1, Zelevinsky showed that the varieties $\line{O}_{\b}$ are locally
isomorphic to some Schubert varieties of type $A_{m}$, where $m=\deg(\b)$. 
Hence we know that $P_{\a, \b}(q)$ 
is a polynomial in $q$. More precisely, for each $\a$, Zelevinsky
associated a permutation $w(\a)$ in the symmetric group $S_{\deg(\a)}$ such that we have 
\[
P_{\b, \a}(q)=P_{w(\a), w(\b)}(q).
\]
Unfortunately, the description of Zelevinsky on the element $w(\a)$ is brutal 
and inexplicite, which lacks geometric meaning.

\remk
In this paper, 
for symmetric multisegments $\a$ and $\b$ (cf. section 4), 
we will give more concrete description about the coefficient $m_{\b, \a}$
in terms of elements in $S_{n}$ with $n$ equals to the number of segments
contained in $\a$, cf. corollary \ref{cor: 2.5.9}.
For general case, we will give use the reduction method from
section \ref{section-reduction-symmetric-case} to give a more concrete description on $P_{\b, \a}(q)$.

\section{Symmetric multisegments and the associated orbital varieties}\label{section-multi-variaties}
In this section we introduce the notion of symmetric 
multisegment, which plays an essential role in our present
paper.
\begin{definition}
Let $\Delta=[i, j]$ be a segment, 
then we define the beginning and the end of $\Delta$ to be
\[
b(\Delta)=i, \quad e(\Delta)=j.
\]
\end{definition}

\begin{definition}
We say a multisegment $\a=\{\Delta_1, \cdots, \Delta_r\}$
is regular if $b(\Delta_1), \cdots, b(\Delta_r)$ are distinct 
and $e(\Delta_1), \cdots, e(\Delta_r)$ are distinct.
\end{definition}

\begin{example}
The segment $\a=\{[1, 2], [2, 4], [4, 5]\}$ is regular.
\end{example}

\begin{prop}
Let $\a$ be a regular multisegment, then any $\b\leq \a$ is 
also regular.
\end{prop}

\begin{proof}
This follows from the fact that if $\a_1$ is obtained 
from $\a$ by elementary operation, then 
$b(\a_1)\subseteq b(\a)$ and $e(\a_1)\subseteq e(\a)$.
\end{proof}

\begin{definition}\label{def: 2.1.5}
 Let $\a=\{\Delta_{1}, \cdots, \Delta_{n}\}$  be regular.
We say that $\a$ is symmetric if 
\[
 \max\{b(\Delta_{i}): i=1,\cdots, n\}\leq \min\{e(\Delta_i): i=1,\cdots, n\}.
\]

\end{definition}

\begin{example}
The multisegment $\a=\{[1, 4], [2, 5], [3,6]\}$ is symmetric.
\end{example}

We have  
\begin{prop}\label{teo: 2.3.2}
Fix a symmetric multisegment $\a_{\Id}=\{\Delta_{1}, \cdots, \Delta_{n}\}$ satisfying 
\[
  b(\Delta_{1})< \cdots < b(\Delta_{n}),
\]

\[
  e(\Delta_{1})< \cdots < e(\Delta_{n}).
\]
Then for any permutation in the symmetric group $S_n$, the formula 
\begin{align*}
\Phi(w)&=\sum_{i=1}^{n} [b(\Delta_{i}), e(\Delta_{w(i)})]
\end{align*}
defines a bijection between $S_{n}$ and $S(\a_{\Id})$.
 Moreover, the order relation on $S(\a_{\Id})$ induces the inverse Bruhat order, i.e.,
 \[
  w\leq v\Leftrightarrow \Phi(w)\geq \Phi(v).
 \]
\end{prop}

\begin{proof}
The injectivity is clear.
We observe that $\Phi(\Id)=\a_{\Id}$. 
We show now that $\Phi(w)\in S(\a_{\Id})$ for general $w$
and the partial order on $S(\a_{\Id})$ induces the inverse Bruhat order.

\begin{enumerate}
\item[(1)]For $v\leq w$,
 by the chain property of Bruhat order(cf. \cite{BF} Theorem 2.2.6), we have
\[
 v=w_{0}<w_{1}<\cdots <w_{\alpha}=w,
\]
such that $w_{\gamma}=\sigma_{i_{\gamma-1}, j_{\gamma-1}}w_{\gamma-1}$ for some
$i_{\gamma-1}< j_{\gamma-1}$
and $\ell(w_{\gamma})=\ell(w_{\gamma-1})+1$. 
Now by lemma 2.1.4 of \cite{BF}, we know that 

\[
 w_{\gamma-1}^{-1}(i_{\gamma-1})<w_{\gamma-1}^{-1}(j_{\gamma-1}).
\]
Hence the segments
$$[b(\Delta_{w_{\gamma-1}^{-1}(i_{\gamma-1})}), e(\Delta_{i_{\gamma-1}})]$$ 
$$[b(\Delta_{w_{\gamma-1}^{-1}(j_{\gamma-1})}), e(\Delta_{j_{\gamma-1}})]$$
are linked in $\Phi(w_{\gamma-1})$. Moreover, by performing the
elementary operation on the two segments, we obtain $\Phi(w_{\gamma})$, 
hence
\[
 \Phi(w_{\gamma-1})>Phi(w_{\gamma}).
\]
Again by transitivity of partial orders, we are done.
Note that we proved that all $\Phi(w)$ are in $S(\a_{\Id})$.
Moreover, for any $\b\in S(\a_{\Id})$, 
the fact that $\a_{\Id}$ is symmetric implies $b(\a_{\Id})=b(\b)$ since no segment 
is juxtaposed to the others. 
The same reason shows that $e(\a_{\Id})=e(\b)$.
Hence 
there is a unique $w\in S_{n}$ 
such that 
\[
\b=\sum_{i=1}^{n}[b(\Delta_{i}), e(\Delta_{w(i)})].
\]
This proves the surjectivity.

\item [(2)] Let $\Phi(w)\geq \Phi(v)$,
 we choose 
$$
 \Phi(w)=\Phi(w_{0})>\Phi(w_{1})>\cdots >\Phi(w_{\alpha})=\Phi(v)
$$
to be a maximal chain of multisegments, where 
$\Phi(w_{\gamma})$ is obtained from $\Phi(w_{\gamma-1})$ by performing the elementary operation 
to segments 
$$
\{[b(\Delta_{i_{\gamma-1}}), e(\Delta_{w_{\gamma-1}(i_{\gamma-1})})], 
\quad [b(\Delta_{j_{\gamma-1}}), e(\Delta_{w_{\gamma-1}(j_{\gamma-1})})]\}
$$
in $\Phi(w_{\gamma-1})$ with $i_{\gamma-1}<j_{\gamma-1}$. Then 
\[
 w_{\gamma-1}(i_{\gamma-1})<w_{\gamma-1}(j_{\gamma-1}).
\]

Hence 
$$
w_{\gamma}=\sigma_{w_{\gamma-1}(i_{\gamma-1}), w_{\gamma-1}(j_{\gamma-1})}w_{\gamma-1}.
$$
Note that in this case, we have either 
\[
 w_{\gamma}<w_{\gamma-1}
\]
or 
\[
 w_{\gamma}>w_{\gamma-1},
\]
by (1), the former implies $\Phi(w_{\gamma-1})<\Phi(w_{\gamma})$, 
contradiction to our assumption.

Hence we must have
\[
 w_{\gamma}>w_{\gamma-1}.
\]
We conclude by transitivity of partial order that $w<v$.
\end{enumerate}
\end{proof}

For the moment, we consider a special case of symmetric multisegments, 
we assume that
\[                       
\a_{\Id}=\sum_{i=1}^{n}[i, n+i-1],                     
\]
with weight function
\[
\quad \varphi=\sum_{\Delta\in \a}f_{\a}(\Delta)\chi_{\Delta}.
\]
We remind that we already constructed a bijection
\[
 \Phi: S_{n}\rightarrow S(\a_{\Id})
\]
such that $\Phi(\Id)=\a_{\Id}$.

We consider the orbital variety $E_{\varphi}$ attached to $\a_{\Id}$.
\begin{definition}
 Let 
 $$O_{w}=O_{\Phi(w)},\quad \hbox{ and }\quad  O_{\varphi}^{\sym}=\coprod_{w\in S_{n}} O_{w}\subseteq E_{\varphi}.$$
 Also, let $$\line{O}_{w}^{sym}=\line{O}_{w}\cap O_{\varphi}^{\sym}.$$
\end{definition}

\begin{definition}
 Let $M_{i, j}$ be the space of $i\times j$ matrices. Let 
 \[
 \xymatrix
 { 
 E_{\varphi}\ar[d]^{ p_{\varphi}}&\hspace{-1cm} =M_{2,1}\times \cdots M_{n-1, n-2}\times M_{n,n-1}\times M_{n-1,n}
 \times \cdots \times M_{1,2}\\
 Z_{\varphi}&\hspace{-2.8cm}: =M_{2,1}\times \cdots M_{n-1, n-2}\times M_{n-1,n}\times \cdots \times M_{1,2}.
 }
 \]
\end{definition}
be the natural projection
with fiber $M_{n,n-1}$.

Now we want to describe the fiber of the 
restriction $p_{\varphi}|_{O_{\varphi}^{\sym}}$.
\begin{definition}
 We define $GL_{n,n-1}$ to be the subset of $M_{n, n-1}$ consisting
of the matrices of rank $n-1$.

\end{definition}

We denote by $p_{n}: M_{n, n}\twoheadrightarrow M_{n, n-1}$ the morphism
of forgetting the last column in $M_{n, n}$.

\remk
Now by restriction to $GL_{n}$, we have the morphism
\[
 p_{n}: GL_{n}\twoheadrightarrow GL_{n,n-1},
\]
which satisfies the property that $p_n(g_1g_2)=g_1p_n(g_2)$ for $g_1, g_2\in GL_n$.

\begin{prop}\label{prop: 4.3.6}
The morphism 
\[
  p_{n}: GL_{n}\twoheadrightarrow GL_{n,n-1},
\]
is a fibration. 
Furthermore, it induces a bijection 
\[
p_{n}: B_{n}\backslash GL_{n}/B_{n}\rightarrow B_{n}\backslash GL_{n,n-1}/B_{n-1}.
\]
\end{prop}

\begin{proof}
To see that it is locally trivial, 
note that $p_n$ is $GL_n$ equivariant with $GL_n$ acting by multiplication
from the left. Since $GL_n$ acts transitively on itself, it acts also 
transitively on $GL_{n, n-1}$. Now $p_n$ is equivariant implies that 
all the fibers of $p_n$ are isomorphic. 
Let $H$ be the stabilizer of $p_n(\Id)$, then 
$GL_{n, n-1}\simeq GL_n/H$, it is a \'etale locally trivial 
fibration according to Serre \cite{S} proposition
 3.
 By Bruhat decomposition, every $g\in GL_{n}$ admits a decomposition
\[
 g=b_{1}wb_{2}, \quad b_{i}\in B_{n}, i=1,2, \quad w\in S_{n},
\]
here we identify $S_n$ with the set of permutation matrices in $GL_n$.
We can decompose $b_{2}=b_{3}v$, where $b_{3}\in GL_{n-1}$, which is identified 
with the direct summand in the Levi subgroup $GL_{n-1}\times \C^{\times}$,
and $v-\Id$ only contains non zero elements
in the last column, 
by definition, 
\[
 p_{n}(g)=b_{1}p_{n}(w)b_{3}.
\]
We obtain that $p_n$ induces 
\[
p_{n}: B_{n}\backslash GL_{n}/B_{n}\rightarrow B_{n}\backslash GL_{n,n-1}/B_{n-1}.
\]

It is bijective
because given $p_{n}(w)$, there is a unique way to reconstruct an
element which belongs to $S_{n}$.
\end{proof}

\begin{teo}\label{teo: 4.3.7}
The morphism
$$
p_{\varphi}|_{O_{\varphi}^{\sym}}
$$
is smooth with fiber $GL_{n,n-1}$. 
Moreover, the morphism $p_{\varphi}|_{O_{w}}: O_{w}\rightarrow p_{\varphi}(O_{\varphi}^{\sym})$ 
is surjective with fiber $B_{n}p_{n}(w)B_{n-1}$.
\end{teo}

\begin{proof}
Note that smoothness follows from that $p_{\varphi}: E_{\varphi}\rightarrow Z_{\varphi}$ 
is smooth and that $O_{\varphi}^{\sym}$ is open in $E_{\varphi}$. 
To see the rest of the properties, we fix an element $e_{w}$ in 
each orbit $O_{w}$ as follow.
Let $(v_{ij})(i=1, \cdots, 2n-1, j=1, \cdots, \varphi(i))$ be a basis of $V_{\varphi, i}$, 
and an element $e_{w}$
satisfying
\begin{displaymath}
    \left\{ \begin{array}{cc}
            e_{w}(v_{ij})&\hspace{-0.5cm}=v_{i+1,j}, \quad\text{ for } i<n-1\\
            e_{w}(v_{n-1, j})&\hspace{-3 cm}=v_{n, w(j)}, \\
            e_{w}(v_{ij})&\hspace{-0.5 cm}=v_{i+1,j-1}, \quad\text{ for } i\geq n.
           \end{array}\right.,
\end{displaymath}
here we put $v_{i,0}=0$. 
\begin{example} 
Let $w=(1,2)$, then by the strategy in the proof, 
$e_{w}$ is given by the following picture:
\begin{figure}[!ht]
\centering
\includegraphics{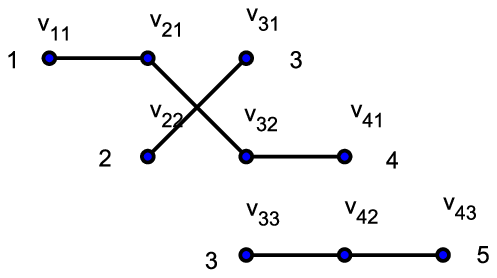}
\caption{\label{fig-multisegment13} }{Construction of $e_{w}$ in case $n=3$}
\end{figure} 
\end{example}
\vspace{0.5cm}
We claim that $e_{w}\in O_{w}$. 
In fact, it suffices to observe that 
\[
 e_{w}: v_{ii}\rightarrow \cdots\rightarrow v_{n-1,i}\rightarrow v_{n, w(i)}\rightarrow v_{n+1, w(i)-1}\rightarrow \cdots 
 v_{n+w(i)-1, 1},
\]
which by proposition
\ref{prop: 2.2.4}, implies that the multisegment indexing $e_{w}$ contains $[i,w(i)+n-1]$
for all $i=1, \cdots, n$, 
hence it must be $\Phi(w)$. 
Note that, by definition, we have 
\[
 p_{\varphi}(e_{\Id})=p_{\varphi}(e_{w}), ~\text{for all  }w\in S_{n}.
\]
Since $p_{\varphi}$ is compatible with the action of $G_{\varphi}$, 
we get 
\[
 p_{\varphi}(O_{\varphi}^{\sym})=p_{\varphi}(O_{w}),~ \text{for all }w\in S_{n},
\]
which implies that $p|_{O_{w}}$ is surjective. Now it remains to characterize
its fiber.
Let $T'\in p_{\varphi}(O_{\varphi}^{\sym})$,  then $p_{\varphi}^{-1}(T')\simeq M_{n, n-1}$
in $E_{\varphi}$. Moreover, for $T=(T_1, \cdots, T_{2n-2})\in p_{\varphi}^{-1}(T')$, then $T\in O_{\varphi}^{\sym}$ 
if and only if 
\[
T_{n-1}\in GL_{n, n-1}.
\]
Therefore, the map $T\mapsto T_{n-1}$ induces 
\[
p_{\varphi}^{-1}(T')\cap O_{\varphi}^{\sym}\simeq GL_{n, n-1}.
\]

Consider the variety $p_{\varphi}^{-1}(T')\cap O_w$. Note that 
since $G_{\varphi}$ acts transitively on  $p_{\varphi}(O_{\varphi}^{\sym})$, we may assume that $T'=p_{\varphi}(e_{\Id})$. 

\begin{lemma}\label{lem: 4.3.8}
The set of $f_{w}\in O_{w}$ satisfying 
\begin{displaymath}
    \left\{ \begin{array}{cc}
            f_{w}(v_{ij})&=v_{i+1,j}, \quad \text{ for } i<n-1\\
            f_{w}(v_{ij})&\hspace{-0.5cm}=v_{i+1,j-1}, \quad \text{ for } i\geq n.
           \end{array}\right.
\end{displaymath}
is in bijection with $B_{n}p_{n}(w)B_{n-1}$ via $p_{\varphi}^{-1}(p_{\varphi}(e_{\Id}))\cap O_{\varphi}^{\sym}\simeq GL_{n, n-1}$.
\end{lemma}

\begin{proof}
%

Now the element $f_{w}\in O_{w}$ is completely determined by 
the component 
\[
 f_{w, n-1}: V_{\varphi, n-1}\rightarrow V_{\varphi, n}.
\]
We know by proposition \ref{prop: 2.2.4} that 
$f_{w, n-1}$ is injective hence of rank $n-1$.
Hence we have 
$f_{w, n-1}\in GL_{n, n-1}$.

Now by proposition \ref{prop: 4.3.6} we get $B_{n}\backslash GL_{n,n-1}/B_{n-1}$
is indexed by $S_{n}$, it remains to see that $f_{w, n-1}$ is in the class indexed by
$p_{n}(w)$. 

Finally, we note that $p_{\varphi}$ is a morphism equivariant 
under the action of 
\[
 G_{\varphi}=GL_{1}\times GL_{2}\times \cdots \times GL_{n-1}\times GL_{n}\times \cdots \times GL_{2}\times GL_{1}.
\]
Since $G_{\varphi}$ acts transitively on $O_{w}$,  the image of $O_{w}$ is $G_{\varphi}. (p_{\varphi}(e_{w}))$, hence is 
$p_{\varphi}(O_{\Id})$. Now we prove that the stabilizer of 
$p_{\varphi}(e_{w})$ is $B_{n}\times B_{n-1}$.
Let $e_{\Id}=(e_{1}, \cdots, e_{ n-1}, e_{ n}, \cdots, e_{2n-2})$
with $e_{i}\in M_{ i, i+1}$ if $i<n$ and $e_i\in M_{i, i-1}$ if $i\geq n$. We have 
$$p_{\varphi}(e_{\Id})=(e_{1}, \cdots, e_{ n-2}, e_{n}, \cdots, e_{2n-2}).$$

Let $g=(g_{1}, \cdots, g_n, g_{n+1}, \cdots, g_{2n-1})$ such that 
$g.p_{\varphi}(e_{\Id})=p_{\varphi}(e_{\Id})$. Then by definition for $i<n-1$
we know that 
$g_{i+1}e_ig_i^{-1}=e_i$. We prove by induction on $i$ that 
$g_i\in B_i\in GL_i$ for $i\leq n-1$. For $i=1$, we have nothing to prove.
Now assume that $i\leq n-2$, and $g_i\in B_i$, we show that $g_{i+1}\in B_{i+1}$. 
Consider $$g_{i+1}e_ig_i^{-1}(g_i(v_{ij}))=g_{i+1}e_i(v_{ij})=g_{i+1}(v_{i+1, j}).$$
On the other hand,
by induction, we know that 
$$g_{i+1}e_ig_i^{-1}(g_i(v_{ij}))=e_i(g_i(v_{ij}))\in \oplus_{k\leq j}\C v_{i+1, k}.$$ 
Therefore we have $g_{i+1}\in B_{i+1}$. Actually, since $e_i$ is injective, 
the equality $e_i(g_i(v_{ij}))=g_{i+1}(v_{i+1, j})$, implies that $g_i$ is 
completely determined by $g_{i+1}$.
This shows that $g_{n-1}\in B_{n-1}$ it determines all $g_i$ for $i\geq n-1$. The same method proves that 
$g_{n}\in B_n$ and it determines all $g_i$ for $i\geq n$. 
We conclude that the fiber of the morphism
$p_{\varphi}|_{O_{w}}$ is isomorphic to 
$B_{n}p_{n}(w)B_{n-1}$.
\end{proof}
\end{proof}

\begin{cor}
We have for $v\leq w$ in $S_{n}$, and $X_{w}$ the closure of $B_{n}wB_{n}$ in $GL_{n}$,
\[
\dim \mathcal{H}^{i}(\line{O}_{w}^{sym})_{v}=\dim\mathcal{H}^{i}(X_{w})_{v},
\]
for all $i\in \Z$, here the index $v$ on the left hand side indicates that we localize at a generic point in $O_{v}$ and 
on the right hand side means that we localize at a generic point in $B_{n}vB_{n}$.
\end{cor}

\begin{proof}
Since $p_{\varphi}|_{O_{\varphi}^{\sym}}$ is a fibration with fiber $GL_{n,n-1}$ 
over $Z_{\varphi}$, we apply the smooth base change theorem to the following Cartesian diagram
\begin{displaymath}
\xymatrix 
{
GL_{n,n-1}\ar[r]\ar[d]&O_{\varphi}^{\sym}\ar[d]\\
p_{\varphi}(\Phi(\Id))\ar[r]&Z_{\varphi}.
}
\end{displaymath}
We get
$$
\dim \mathcal{H}^{i}(\line{O}_{w}^{sym})_{v}=\dim \mathcal{H}^{i}(\line{B_np_n(w)B_{n-1}})_{B_np_n(v)B_{n-1}}.
$$ 
Now apply proposition \ref{prop: 4.3.6}, we have 
\[
\dim \mathcal{H}^{i}(\line{B_np_n(w)B_{n-1}})_{B_np_n(v)B_{n-1}}=dim\mathcal{H}^{i}(X_{w})_{v}.
\]
\end{proof}

\begin{cor}\label{cor: 2.5.9}
 We have for $v\leq w$ in $S_{n}$,
 \[
  m_{\Phi(v), \Phi(w)}=P_{v, w}(1).
 \]
 \end{cor}
 
\begin{proof}
 This follows from the fact that 
 \[
  \dim \sum_i\mathcal{H}^{i}(X(w))_{v}=P_{v,w}(1)
 \]
(cf. \cite{KL}).
\end{proof}

\section{Descent of Degrees for Multisegment}

In this section we describe a procedure to 
decrease the degree of a multisegment $\a$ without
affecting the coefficients $m(\b, \a)$.

\subsection{Notation and Combinatorics}

\begin{notation}\label{nota: 3.1.2}
 For $\Delta=[i,j]$ a segment, we put 
 \begin{align*}
 \Delta^{-}=&[i,j-1], \quad {^{-}\Delta}=[i+1,j],\\
 \Delta^{+}=&[i,j+1], \quad {^{+}\Delta}=[i-1,j].
 \end{align*}
\end{notation}

\begin{definition}\label{def: 2.2.3}
Let $k\in \Z$ and $\Delta$ be a segment, we define 
\begin{displaymath}
 \Delta^{(k)}=\left\{ \begin{array}{cc}
         &\Delta^-, \text{ if } e(\Delta)=k;\\
         &\hspace{-0.5cm}\Delta, \text{ otherwise }.
         \end{array}\right.
\end{displaymath}
For a multisegment 
 $
 \a=\{\Delta_{1}, \cdots, \Delta_{r}\},
 $
we define
  \[
  \a^{(k)}=\{\Delta_{1}^{(k)}, \cdots,\Delta_{r}^{(k)} \}.
 \]
\end{definition}

\begin{definition}\label{def: 3.1.3}
We say that the multisegment $\b\in S(\a)$ satisfies\textbf{ the hypothesis $H_{k}(\a)$} if
the following two conditions are verified 
\begin{description}
 \item[(1)]$\deg(\b^{(k)})=\deg(\a^{(k)})$;
 \item[(2)]there exists no pair of linked segments $\{\Delta, \Delta'\}$ in $\b$ such that $e(\Delta)=k-1, ~e(\Delta')=k$.
\end{description}
\end{definition}

\begin{definition}
 Let 
\[
 \tilde{S}(\a)_{k}=\{\c\in S(\a): \deg(\c^{(k)})=\deg(\a^{(k)})\}.
\]
\end{definition}

\begin{lemma}\label{lem: 3.1.5}
Let $\c\in \tilde{S}(\a)_{k}$. Then  
\[
 \sharp\{\Delta\in \a: e(\Delta)=k\}=\sharp\{\Delta\in \c: e(\Delta)=k\}.
\] 
\end{lemma}

\remk
Here we count segments with multiplicities.

\begin{proof}
Note that 
\[
 \deg(\a)=\deg(\a^{(k)})+\sharp\{\Delta\in \a: e(\Delta)=k\}.
\]
Now that for $\c\in \tilde{S}(\a)_{k}$
\[
 \deg(\c)=\deg(\a),\quad \deg(\c^{(k)})=\deg(\a^{(k)}),
\]
we have 
\[
 \sharp\{\Delta\in \a: e(\Delta)=k\}=\sharp\{\Delta\in \c: e(\Delta)=k\}.
\]
\end{proof}

 \begin{lemma}\label{lem: 3.0.7}
 Let $k\in \Z$.
 \begin{description}
  \item[(1)] For any $\b\in S(\a)$, we have $\deg(\b^{(k)})\geq \deg(\a^{(k)})$.
 \item[(2)]Let $\c\in \tilde{S}( \a)_{k}$, then for $\b\in S(\a)$ such that 
 $\b>\c$, we have $\b\in \tilde{S}(\a)_{k}$.
 \item[(3)]Let $\b\in \tilde{S}(\a)_k$, then $\b^{(k)}\in S(\a^{(k)})$. Moreover, if we suppose that $\a$ satisfies the hypothesis $H_{k}(\a)$
  and $\b\neq \a$, then 
  $$\b^{(k)}\in S(\a^{(k)})-\{\a^{(k)}\}$$.
  
 \item[(4)]Suppose that $\a$ does not verify the hypothesis $H_{k}(\a)$, then 
 there exists $\b\in S(\a)$ satisfying the hypothesis $H_{k}(\a)$, such that
 $\b^{(k)}= \a^{(k)}$.
 \end{description}
 \end{lemma}

\begin{proof}
For (1), note that 
for any $\b\in S(\a)$, $e(\b):=\{e(\Delta): \Delta\in \b\}$ is a sub-multisegment
of $e(\a)$. And from $\b$ to $\b^{(k)}$, 
we replace those segments $\Delta$ such that $e(\Delta)=k$
by $\Delta^{-}$. Now (1) follows by counting the segments ending 
in $k$.

For (2), by (1), we have 
\[
 \deg (\a^{(k)})\leq \deg(\b^{(k)})\leq \deg(\c^{(k)}).
\]
The fact that $\c\in \tilde{S}( \a)_{k}$ implies that 
$ \deg (\a^{(k)})= \deg (\c^{(k)})$, hence 
$ \deg (\a^{(k)})= \deg (\b^{(k)})$ and $\b\in \tilde{S}( \a)_{k}$.

As for (3), 
suppose that 
$\deg(\b^{(k)})=\deg(\a^{(k)})$, we prove 
$\b^{(k)}<\a^{(k)}$.
Let
\[
 \a=\a_{0}>\cdots>\a_{r}=\b
\]
be a maximal chain of multisegments, then 
by $(2)$, we know 
$\deg(\a_{j}^{(k)})=\deg(\a^{(k)})$, for all $ j=1,\cdots, r$.
Our proof breaks into two parts.

\vspace{0.5cm}
(I)We show that 
\[
 \deg(\a_{j}^{(k)})=\deg(\a_{j+1}^{(k)})\Rightarrow \a_{j}^{(k)}\geq \a_{j+1}^{(k)}.
\]

Let $\a_{j+1}$ be obtained from $\a_{j}$
by applying the elementary operation to two linked segments $\Delta, \Delta'$.
\begin{itemize}
\item 
If none of them ends in $k$, then $\a_{j}^{(k)}$ contains
both of them. And we obtain $\a_{j+1}^{(k)}$ by applying 
the elementary operation to them.

If one of them ends in $k$, 
we assume $e(\Delta')=k$.

\item If 
$\Delta$ precedes $\Delta'$,  we know that if $e(\Delta)<k-1$,
$\Delta$ is still linked to $\Delta'^-$, and one obtains $\a_{j+1}^{(k)}$
by applying elementary operation to $\{\Delta, ~\Delta'^-\}$,
 otherwise $e(\Delta)=k-1$, which implies 
 $\a_{j+1}^{(k)}=\a_{j}^{(k)}$.

\item If $\Delta$ is preceded by $\Delta'$, 
then the fact that 
\[
 \deg(\a_{j+1}^{(k)})=\deg(\a_{j}^{(k)})
\]
implies $b(\Delta)\leq k$, hence
$\Delta'^{-}$ is linked to $\Delta$, and we obtain
$\a_{j+1}^{(k)}$ from $\a_{j}^{(k)}$ by applying elementary
operation to them.
\end{itemize}
Here we conclude that $\b^{(k)}\in S(\a^{(k)})$.
\vspace{0.5cm}

(II)Assuming that $\a$ satisfies the hypothesis $H_k(\a)$, we show that 
\[
 \a_{1}^{(k)}<\a^{(k)}.
\]
Let $\a_{1}$ be obtained from $\a$
by performing the elementary operation to $\Delta, \Delta'$.

We do it as in (1) but put $j=0$. Note that in (1), the only case where we can have 
$\a_{1}^{(k)}=\a^{(k)}$ is when 
$\Delta$ precedes $\Delta'$  and $ e(\Delta')=k, ~e(\Delta)=k-1$.
But such a case can not exist  since $\a$ verifies the hypothesis $H_{k}(\a)$.
Hence we are done.

Finally, for (4), we construct $\b$ in the following 
way. Suppose that $\a$ does not satisfy the 
hypothesis $H_{k}(\a)$, then there exists a pair of 
linked segments $\{\Delta, \Delta'\}$ such that 
\[
 e(\Delta)=k-1, \quad e(\Delta')=k,
\]
let $\a_{1}$ be the multisegment obtained by applying the
elementary operation to $\Delta$ and $\Delta'$. We have
\[
 \a_{1}^{(k)}=\a^{(k)}.
\]
If again $\a_{1}$ fails the hypothesis $H_{k}(\a)$, we repeat the same
construction to get $\a_{2}, \cdots$, since
\[
 \a>\a_{1}>\cdots. 
\]
In finite step, we get $\b$ satisfying the conditions in the theorem and 
\[
 \b^{(k)}=\a^{(k)}.
\]
\end{proof}
\remk Actually, the multisegment constructed in $(4)$ is unique, as we shall see 
later(proposition \ref{cor: 3.2.3}).

\begin{definition}
We define a morphism 
\[
 \psi_{k}:\tilde{S}( \a)_{k}\rightarrow S(\a^{(k)})
\]
by sending $\c$ to $\c^{(k)}$. 
\end{definition}

\begin{prop}\label{prop: 1.5.8}
The morphism $\psi_{k}$ is surjective. 
 
\end{prop}

\begin{proof}
Let $\d\in S(\a^{(k)})$, such that 
we have a maximal chain of multisegments, 
\[
 \a^{(k)}=\d_{0}>\cdots >\d_{r}=\d.
\]
By induction, we can assume that 
there exists $\c_{i}\in \tilde{S}(\a)_{k}$
such that $\c_{i}^{(k)}=\d_{i}$, for all 
$i<r$. Assume we obtain $\d$ from $\d_{r-1}$ by performing the 
elementary operation on the pair of linked segments $\{\Delta\prec \Delta'\}$.

\begin{itemize}
 \item If $e(\Delta)\neq k-1$ and $e(\Delta')\neq k-1$, then 
 we observe that the pair of segments are actually contained in $\c_{r-1}$.
 Let $\c_{r}$ be the multisegment obtained by performing the elementary operation 
 to them . We conclude that $\c_{r}^{(k)}=\d_{r}$, and $\c\in \tilde{S}(\a)_{k}$.
 \item If $e(\Delta)=k-1$, then $\Delta\in \c_{r-1}$ or $\Delta^+\in \c_{r-1}$ and $\Delta'\in \c_{r-1}$.
 The fact that $\d_{r-1}=\c_{r-1}^{(k)}$
implies that $k\notin e(\d_{r-1})$, hence $e(\Delta')>k$. Hence 
both $\Delta$ and $\Delta^{+}$ are linked to $\Delta'$. In either case we perform 
the elementary operation to get $\c_{r}$ such that $\c_{r}^{(k)}=\d$. 
\item If $e(\Delta')=k-1$, then $\Delta'\in \c_{r-1}$ or $\Delta'^+\in \c_{r-1}$ and $\Delta\in \c_{r-1}$.
The same argument as in the second case shows that there exists $\c_{r}$ such that 
$\c_{r}^{(k)}=\d$.
\end{itemize}
 
\end{proof}

Actually, the proof in proposition \ref{prop: 1.5.8} 
yields the following refinement.
\begin{cor}\label{cor: 1.5.9}
Let $\c\in \tilde{S}(\a)_{k}, \d\in S(\a^{(k)})$ such that 
\[
 \c^{(k)}>\d,
\]
then there exists a multisegment $\e\in\tilde{S}(\a)_{k}$ such that
\[
 \c>\e, ~ \e^{(k)}=\d.
\]
\end{cor}

\begin{proof}
Note that $\c\in \tilde{S}(\a)_k$ implies $\tilde{S}(\a)_k\supseteq \tilde{S}(\c)_k$.
Combine with the surjectivity of 
\[
\psi_k: \tilde{S}(\c)_k\rightarrow S(\c^{(k)}),
\]
we get the result.
\end{proof}

\begin{definition}\label{def: 2.3.2}
For $\a$ a multisegment, and $k\in \Z$
we define 
\[
 S(\a)_{k}=\{\c\in \tilde{S}(\a)_{k}:
 \c \text{ satisfies the hypothesis } H_{k}(\a) \}.
\]
\end{definition}

%
%
%

%

\begin{prop}\label{prop: 1.5.10}
 The restriction 
 \begin{align*}
  \psi_{k}: S(\a)_{k}&\rightarrow S(\a^{(k)})\\
  \c&\mapsto \c^{(k)} 
 \end{align*}
is also surjective.
\end{prop}

\begin{proof}
 For $\d\in S(\a^{(k)})$, by proposition \ref{prop: 1.5.8},
 we know that there exists $\c\in \tilde{S}(\a)_{k}$
 such that $\c^{(k)}=\d$. 
 But by $(4)$ in lemma \ref{lem: 3.0.7}, we know that 
 there exists $\c'\in S(\c)_{k}$ such that 
 $\c'^{(k)}=\c^{(k)}=\d$. 
 We conclude by the observation that if $\c\in \tilde{S}(\a)_{k}$, then 
 \[
  S(\c)_{k}\subseteq S(\a)_{k}.
 \] 
\end{proof}

Also, concerning the corollary \ref{cor: 1.5.9}, we have the 
following 

\begin{cor}\label{cor: 3.1.11}
Let $\c\in \tilde{S}(\a)_{k}$ and $\d\in S(\a^{(k)})$ such that
$\c^{(k)}>\d$. Then there exists a multisegment $\e\in S(\c)_{k}$
such that $\e^{(k)}=\d$.
\end{cor}

\begin{proof}
By corollary \ref{cor: 1.5.9}, we know that there exists 
an $\e'\in \tilde{S}(\c)_{k}$ such that $\e'^{(k)}=\d$.
By $(4)$ in lemma \ref{lem: 3.0.7}, we know that 
 there exists $\e\in S(\e')_{k}$ such that 
 $\e^{(k)}=\e'^{(k)}=\d$. 
 Hence we conclude by 
 the fact that
 if $\e'\in \tilde{S}(\a)_{k}$, then 
 \[
  S(\e')_{k}\subseteq S(\a)_{k}.
 \]

\end{proof}

\begin{definition}
Let $k\in \Z$ and $\Delta$ be a segment.
\begin{displaymath}
 {^{(k)}\Delta}=\left\{ \begin{array}{cc}
         &{^{-}\Delta}, \text{ if } b(\Delta)=k;\\
         &\Delta, \text{ otherwise }.
         \end{array}\right.
\end{displaymath}

Let
$$
\a=\{\Delta_{1}, \cdots, \Delta_{r}\},
$$
be a multisegment, we define 
  \[
  {^{(k)}\a}=\{{^{(k)}\Delta_{1}}, \cdots,{^{(k)}\Delta_{r}}, \}.
 \]
\end{definition}

\begin{definition}  We say that the multisegment $\b\in S(\a)$ satisfies \textbf{ the hypothesis ${_{k}H(\a)}$} if
 the following two conditions are verified 
\begin{description}
 \item[(1)]$\deg({^{(k)}\b})=\deg({^{(k)}\a})$;
 \item[(2)] there exists no pair of linked segments $\{\Delta, \Delta'\}$ such that 
 $$b(\Delta)=k, ~b(\Delta')=k+1.$$
\end{description}
\end{definition}

\remk 
There exists a version of lemma \ref{lem: 3.0.7} for ${^{(k)}a}$.
In the following sections, we will work exclusively with 
$\a^{(k)}$ and the hypothesis  $H_{k}(\a)$.
But all our results will remain valid if we replace $\a^{(k)}$ 
by  ${^{(k)}\a}$ and $H_{k}(\a)$ by ${_{k}H(\a)}$.

\subsection{Injectivity of \texorpdfstring{$\psi_{k}$}{Lg}: First Step}
By previous section, we know there exists $\c\in S(\a)_k$ such that $\c^{(k)}=(\a^{(k)})_{\min}$, the 
minimal element in $S(\a^{(k)})$. In this section, we give an explicit construction of such a $\c$
and show that it is the unique multisegment in $S(\a)_k$ which is sent to $(\a^{(k)})_{\min}$ by $\psi_k$.

\begin{itemize}
\item[(i)] In proposition \ref{prop: 3.2.2}, we construct a multisegment $\c\in S(\a_1)_k$ such that $\c^{(k)}=(\a^{(k)})_{\min}$, 
where $\a_1$ is a multisegment such that $\a\in S(\a_1)$.

\item[(ii)] We prove that if there exists a multisegment in  $S(\a)_{k}$ which is sent to $(\a^{(k)})_{\min}$ by $\psi_k$, then it is unique. 

\item[(iii)] Then we apply the uniqueness result to $S(\a_{1})_{k}$ to prove that the $\c$ in (i) belongs to $S(\a)_k$.\footnote{
Here we use partial derivative to prove our result, but it can also be done 
in a purely combinatorial way, which is less elegant and more lengthy though.
}
\end{itemize}

\begin{notation}\label{nota: 3.2.1}
Let $\ell_{\a, k}=\sharp\{\Delta\in \a: e(\Delta)=k\}$.
\end{notation}


\begin{definition}
Let 
\[
 \a_{0}=\{\Delta\in (\a^{(k)})_{\min}: e(\Delta)=k-1\}.
\]
\end{definition}

\begin{prop}\label{prop: 3.2.2}
Let $\a_{0}=\{\Delta_{1}\succeq \cdots \succeq \Delta_{r}\}$. 
Let $\c$ be a multisegment such that 
\begin{description}
 \item[(1)]If $\varphi_{\a}(k-1)> \varphi_{\a}(k)$, then $r=\varphi_{\a}(k-1)-\varphi_{\a}(k)+\ell_{\a, k}$.
 Let 
 \[
  \c=((\a^{(k)})_{\min}\setminus \a_{0})\cup\{\Delta_{1}^+\succeq \cdots \succeq \Delta_{\ell_{\a, k}}^+\succeq \Delta_{m+1}\succeq \cdots \succeq \Delta_{r} \}.
 \]
\item[(2)]If $\varphi_{\a}(k)-\ell_{\a, k}<\varphi_{\a}(k-1)\leq \varphi_{\a}(k)$, then $r=\varphi_{\a}(k-1)-\varphi_{\a}(k)+\ell_{\a, k}$. Let 
\[
 \c=((\a^{(k)})_{\min}\setminus \a_{0})\cup \{\Delta_{1}^+\succeq \cdots \succeq \Delta_{r}^+\succ\underbrace{ [k]=\cdots =[k]}_{\ell_{k,\a}-r}\}
\]
\item[(3)]If $\varphi_{\a}(k-1)\leq \varphi_{\a}(k)-\ell_{\a, k}$, then $\a_{0}=\emptyset$ and
\[
 \c=\a^{(k)}+\ell_{\a, k}[k].
\]
Then $\c$ satisfies the hypothesis $ H_{k}(\c)$ and $\c^{(k)}=(\a^{(k)})_{\min}$.
\end{description}
\end{prop}

%

\begin{proof}
We show only the case $\varphi_{\a}(k-1)>\varphi_{\a}(k)$,
the proof for other cases is similar.
Note that we have the following equality
\[
 \varphi_{\a}(k-1)=\varphi_{(\a^{(k)})_{\min}}(k-1)=r+\sharp \{\Delta\in (\a^{(k)})_{\min}: \Delta\supseteq [k-1, k]\}.
\]
Moreover, $\varphi_{\a}(k-1)>\varphi_{\a}(k)$ implies 
that no segment in $(\a^{(k)})_{\min}$ starts at $k$ by minimality, hence we also have
\[
 \varphi_{\a}(k)=\varphi_{(\a^{(k)})_{\min}}(k)+\ell_{k,\a}=\sharp \{\Delta\in (\a^{(k)})_{\min}: \Delta\supseteq [k-1, k]\}+\ell_{k,\a}.
\]
Now comparing the two formulas gives
the equality 
$r=\varphi_{\a}(k-1)-\varphi_{\a}(k)+\ell_{\a, k}$.
By definition we have $\c^{(k)}=(\a^{(k)})_{\min}$.
To check that $\c$ satisfies the hypothesis $H_{k}(\c)$, it suffices 
to note that $(\a^{(k)})_{\min}\setminus \a_{0}$ does not contain segment 
which ends in $k-1$.
\end{proof}

\begin{lemma}
 Let $\c\in S(\c)_{k}$ be a multisegment such that $\c^{(k)}$ is minimal. Then 
 if $\d\in S(\c)$ such that $\d^{(k)}=\c^{(k)}$, then $\c=\d$ 
\end{lemma}

\begin{proof}
Suppose that $\d<\c$ is a multisegment such that $\d^{(k)}=\c^{(k)}$. Consider the maximal chain of multisegments
\[
\c=\c_0>\cdots > \c_t=\d.
\]
Our assumption implies that $\c_i^{(k)}=\c^{(k)}$ for all $i=1, \cdots, t$ by lemma \ref{lem: 3.0.7}. Hence we can assume $t=1$ and consider 
$\d\in S(\c)$ to be a multisegment obtained by applying the elementary
operation to the pair of linked segments $\{\Delta\prec \Delta'\}$.
\begin{itemize}
 \item If $e(\Delta)\neq k, e(\Delta')\neq k$, then 
the pair $\{\Delta, \Delta'\}$ also appears in $\c^{(k)}$, 
contradicting the fact that $\c^{(k)}$ is minimal.
\item If $e(\Delta')=k$, then by the fact that $\c\in S(\c)_{k}$, 
we know that $e(\Delta)<k-1$, which implies that the pair 
$\{\Delta, \Delta'^-\}$ is linked and belongs to $\c^{(k)}$,contradiction.
\item If $e(\Delta)=k$ and $b(\Delta')<k+1$, then the 
pair $\{\Delta^-, \Delta'\}$ is still linked and belongs to 
$\c^{(k)}$, contradiction.
\end{itemize} 
Hence we must have $e(\Delta)=k$ and $b(\Delta')=k+1$, this
implies that $\deg(\d^{(k)})>\deg(\c^{(k)})$ and $\d\notin \tilde{S}(\c)_{k}$.
Finally, (b) of lemma \ref{lem: 3.0.7} implies that for all $\d<\c$, 
we have $\d\notin \tilde{S}(\c)_{k}$.
\end{proof}

\begin{prop}\label{prop: 3.2.4}
Let $\c\in S(\c)_{k}$ be a multisegment such that $\c^{(k)}$ is minimal. Then
the partial derivative $\D^{k}(L_{\c})$ contains in $\mathcal{R}$ a unique term 
of minimal degree $L_{\c^{(k)}}$, which appears with multiplicity one.
\end{prop}

\begin{proof}
Let 
$\c=\{\Delta_{1}, \cdots, \Delta_{r}\}$
such that $e(\Delta_{t})=k$ if and only if 
$t=i, \cdots, j$ with $i\leq j$.
Then 
\[
 \D^k(\pi(\c))=\Delta_{1}\times \cdots 
 \times \Delta_{i-1}\times (\Delta_{i}+\Delta_{i}^-)\times
 \cdots \times (\Delta_{j}+\Delta_{j}^-)\times \Delta_{j+1}\times 
 \cdots\times \Delta_{r}
\]
with minimal degree term given by 
\[
 \pi(\c^{(k)})=\Delta_{1}\times \cdots \times \Delta_{i-1}\times
 \Delta_{i}^-\times \cdots \times \Delta_{j}^-\times \Delta_{j+1}
 \times \cdots \times \Delta_{r}.
\]
The same calculation shows that for any $\d\in S(\c)$, 
the minimal degree term in $\D^k(\pi(\d))$ is given 
by $\pi(\d^{(k)})$, whose degree is strictly greater
than that of $\c^{(k)}$ since by previous lemma we know that $\d\notin \tilde{S}(\c)_k$. Note that 
$\D^k(L_{\d})$ is a non-negative sum of irreducible representations
( Theorem \ref{teo: 3}), which do not contain 
any representation of degree equal to that of $\c^{(k)}$, 
by comparing the minimal degree term in 
$\D^k(\pi(\d))$ and $\sum_{\e\in S(\d)}m(\e, \d)\D^k(L_{\e})$.
Finally, comparing the minimal degree term in 
 $\D^k(\pi(\c))$ and $\sum_{\e\in S(\c)}m(\e, \c)\D^k(L_{\e})$
 gives the proposition.
 
\end{proof}

\begin{prop}\label{prop: 3.2.5}
Let $\a$ be a multisegment. 
Then  $S(\a)_{k}$ contains a unique 
multisegment $\c$ such that $\c^{(k)}=(\a^{(k)})_{\min}$.
\end{prop}

\begin{proof}
Let 
$\a=\{\Delta_{1}, \cdots, \Delta_{s}\}$
such that $e(\Delta_{t})=k$ if and only if 
$n=i, \cdots, j$ with $i\leq j$.
Then 
\[
 \D^k(\pi(\a))=\Delta_{1}\times \cdots 
 \times \Delta_{i-1}\times (\Delta_{i}+{\Delta_{i}}^{-})\times
 \cdots \times (\Delta_{j}+\Delta_{j}^-)\times \Delta_{j+1}\times 
 \cdots\times \Delta_{s}
\]
with minimal degree term given by 
\[
 \pi(\a^{(k)})=\Delta_{1}\times \cdots \times \Delta_{i-1}\times
 \Delta_{i}^-\times \cdots \times \Delta_{j}^-\times \Delta_{j+1}
 \times \cdots \times \Delta_{r}.
\]
Note that in $\pi(\a^{(k)})$, 
$m((\a^{(k)})_{\min}, \a^{(k)})=1$(cf. \cite{Z3}). Now compare 
with the terms of minimal degree in 
$\sum_{\d\in S(\a)}m(\d, \a)\D^k(L_{\d})$ 
and apply the proposition \ref{prop: 3.2.5}
yields the uniqueness of $\c$ such that 
$\c^{(k)}=(\a^{(k)})_{\min}$.
\end{proof}

\begin{prop}
 Let $\c$ be the multisegment constructed in 
 proposition \ref{prop: 3.2.2}. Then $\c\in S(\a)$.
\end{prop}

\begin{proof}
Let $$\a_{1}=\a^{(k)}+m[k],$$ then we observe that 
$\a\in S(\a_{1})$.
Because of 
$\c\in S((\a^{(k)})_{\min}+m[k])$, 
we have $\c\in S(\a_{1})$.
Note that since $\deg ((\a_{1})^{(k)})=\deg (\c^{(k)})$,
the fact that $\c\in S(\c)_{k}$ implies that
$\c\in S(\a_{1})_{k}$. Now 
let $\d\in S(\a)_{k}$, then we have 
$\d\in S(\a_{1})_{k}$ since $\deg(\d^{(k)})=\deg (\a_{1}^{(k)})=\deg(\a^{(k)})$.
Assume furthermore that $\d^{(k)}$ is minimal, then 
by  proposition \ref{prop: 3.2.5}, we know that such a multisegment
in $S(\a_{1})_{k}$ is unique, which implies $\d=\c$.
\end{proof}

\begin{cor}\label{cor: 3.2.7}
Let $\c\in S(\a)_{k}$ such that $\c^{(k)}=(\a^{(k)})_{\min}$, then 
$\c$ is minimal in $\tilde{S}(\a)_{k}$. 
\end{cor}

\begin{proof}
 By corollary \ref{cor: 3.1.11}, we know that for any $\d\in \tilde{S}(\a)_{k}$,
 there exists a multisegment $\c'\in S(\a)_{k}$ 
 with $\c'^{(k)}=(\a^{(k)})_{\min}$, such that 
 $\d>\c'$. By uniqueness, we must have $\c=\c'$.
\end{proof}

\subsection{Geometry of Nilpotent Orbits: General Cases}

In this section, we show geometrically that 
the morphism 
\begin{align*}
 \psi_{k}: S(\a)_{k}&\rightarrow S(\a^{(k)})\\
 \c&\mapsto \c^{(k)}  
\end{align*}
is bijective, satisfying the properties
\begin{description}
 \item[(1)]For $\c\in S(\a)_{k}$, we have 
 $m(\c, \a)=m(\c^{(k)}, \a^{(k)})$.
 \item [(2)] The morphism $\psi_{k}$ preserves the order, i.e, 
 for $\c, \d\in S(\a)_{k}$, $\c>\d$ if and only if $\c^{(k)}>\d^{(k)}$.
\end{description}

To achieve this, firstly we consider the sub-variety
$X_{\a}^{k}= \coprod_{\c\in \tilde{S}( \a)_{k}}O_{\c}$, 
and construct a fibration $\alpha$ from $X_{\a}^{k}$
to $Gr(\ell_{\a, k}, V_{\varphi_{\a}, k})$, the latter is the space of the 
$\ell_{\a, k}$-dimensional subspace of $V_{\varphi_{\a}, k}$. 
Secondly, we construct an open immersion 
\[
 \tau_{W}: (X_{\a}^{k})_{W}\rightarrow Y_{\a^{(k)}}\times \Hom(V_{\varphi_{\a}, k-1}, W),
\]
where $(X_{\a}^{k})_{W}$ is the fiber over $W$ with respect to $\alpha$
and $Y_{\a^{(k)}}=\coprod_{\c\in S(\a^{(k)})}O_{\c}$.
Our main difficulty here lies in proving that $\tau_{W}$ is actually an 
open immersion. The idea is to  apply Zariski Main theorem, to do this, 
we have to prove the normality and irreducibility of both varieties. 
Irreducibility of $(X_{\a}^{k})_{W}$ follows from our results 
in previous section, and normality follows from the fibration $\alpha$ 
and the fact that orbital varieties are locally isomorphic to some 
Schubert varieties, by Zelevinsky, cf. \cite{Z4}.

Once we prove that $\tau_{W}$ is an open immersion.
All the desired properties of $\psi_{k}$ then follow.

\vspace{0.5cm}
Here we fix a multisegment $\a$ and let $\varphi=\varphi_{\a}$.

\begin{definition}\label{def: 3.3.1}
\begin{itemize}
 \item 
Let 
$$
X_{\a}^{k}= \coprod_{\c\in \tilde{S}( \a)_{k}}O_{\c},  
$$ 
\item Let
$Y_{\a^{(k)}}= \coprod_{c\in S(\a^{(k)})}O_{\c}$.

\item For $\b>\c$ in $\tilde{S}(\a)_{k}$,
we define 
\[
 X_{\b, \c}^{k}=\coprod_{\b\geq \d\geq\c}O_{\d}.
\] 
\end{itemize}
\end{definition}

Let $\c\in \tilde{S}( \a)_{k}, T\in O_{\c}$, then 
 
\begin{lemma}
Let $\varphi=\varphi_{\a}$.
We have 
$\dim(\ker(T|_{V_{\varphi, k}}))=\sharp\{\Delta\in \a: e(\Delta)=k\}=\ell_{\a, k}$(Notation \ref{nota: 3.2.1}),
which does not depend on the choice of $T$.
\end{lemma}

\begin{proof}
The fact $T\in O_{\c}$ implies 
\[
 \dim(\ker(T|_{V_{\varphi, k}}))=\sharp\{\Delta\in \c: e(\Delta)=k\}.
\]
Then our lemma follows from lemma \ref{lem: 3.1.5}.

\end{proof}

\begin{definition}
Let 
\[
 Gr(\ell_{\a, k}, V_{\varphi})=\{W\subseteq V_{\varphi, k}: \dim(W)=\ell_{\a, k}\},
\]
and for $W\in  Gr(\ell_{\a, k}, V_{\varphi})$, let 
\[
 V_{\varphi}/W=V_{\varphi, 1}\oplus \cdots V_{\varphi, k-1}
 \oplus V_{\varphi, k}/W \oplus \cdots.
\]
Also, we denote by 
 \[
 p_{W}: V_{\varphi}\rightarrow V_{\varphi}/W
 \]
the canonical projection.
\end{definition}

\begin{definition}
We define
\[
 \tilde{Z}^{k}=\{(T, W): W\in Gr(\ell_{\a, k}, V_{\varphi}), T\in End(V/W) \text{ of degree +1}\},
\]
and the canonical projection
\begin{align*}
 \pi:&\tilde{Z}^{k}\rightarrow Gr(\ell_{\a, k}, V_{\varphi})\\
 &(T, W)\mapsto W.
\end{align*}

\end{definition}

\begin{prop}
 The morphism $\pi$ is a fibration with fiber 
 $$
 E_{\varphi_{\a^{(k)}}}.
$$
 \end{prop}

\begin{proof}
 This follows from the definition.
\end{proof}

\begin{definition}\label{def: 3.3.6}
Assume $\b, \c\in S(\a^{(k)})$.
\begin{itemize}
 \item Let
\[
 Z^{k, \a}=\{(T,W)\in \tilde{Z}^{k}: T\in Y_{\a^{(k)}} \}.
\]
\item Let
\[
 Z^{k, \a}_{\b, \c}=\{(T,W)\in \tilde{Z}^{k}: T\in \coprod_{\b\geq\d\geq \c}O_{\d}\}, 
 ~Z^{k, \a}_{\b}=\{(T,W)\in \tilde{Z}^{k}: T\in \coprod_{\d\geq \b}O_{\d}\}.
\]

\item Let
\[
 Z^{k, \a}(\c)=\{(T, W)\in Z^{k, \a}, T\in O_{\c}\}.
\]

\end{itemize}
\end{definition}

\remk 
 The restriction of $\pi$ to $Z^{k, \a}$ is a fibration 
 with fiber $ Y_{\a^{(k)}}$.

\begin{definition}\label{def: 3.3.7}
Now we define $T^{(k)}\in End(V/\ker(T|_{V_{\varphi, k}}))$ such that 
\begin{displaymath}
 T^{(k)}|_{V_{\varphi,i}}=\left\{ \begin{array}{cc}
&\hspace{-1cm}T|_{V_{\varphi, i}}, \text{ for }i \neq k, k-1,\\
&\hspace{-0.4cm}p_{T, k}\circ T|_{V_{\varphi, i}}, \text{ for }i=k-1\\
&\hspace{-1cm} T|_{V_{\varphi, i}}\circ p_{T, k}, \text{ for }i=k.
\end{array}
\right.
\end{displaymath}
where $p_{T, k}: V_{\varphi}\rightarrow V_{\varphi}/\ker(T|_{V_{\varphi, k}})$ is the canonical projection.
\end{definition}

This gives naturally an element $(T^{(k)}, \ker(T|_{V_{\varphi, k}}))$ in $Z^{k, \a}$.
We construct a morphism 
\[
 \gamma_{k}: X_{\a}^{k}\rightarrow Z^{k, \a}.
\]

by 
\[
 \gamma_{k}(T)=(T^{(k)}, \ker(T|_{V_{\varphi, k}})).
\]

\begin{definition}
 We define 
\[
 \alpha: X_{\a}^{k}\rightarrow Gr(\ell_{\a, k}, V_{\varphi}), 
\]
with $\alpha(T)=\ker(T|_{V_{\varphi, k}})$.
\end{definition}

\remk 
 We have a commutative diagram
 \begin{displaymath}
  \xymatrix{
  X_{\a}^{k}\ar[d]^{ \alpha}\ar[r]^{\gamma_{k}}&Z^{k, \a}\ar[dl]^{\pi}\\
  Gr(\ell_{\a, k}, V_{\varphi}).
  }
  \end{displaymath}
where $\gamma_{k}$ maps fibers to fibers.

\begin{prop}\label{prop: 4.6.10}
 The morphism $\alpha$ is a fiber bundle such that $\alpha|_{O_{\c}}$
 is surjective for any $\c\in \tilde{S}(\a)_{k}$.
\end{prop}

\begin{proof}
 We have to show that $\alpha$ is locally trivial. 
 We fix $W\in Gr(\ell_{\a, k}, V_{\varphi})$
 Note that $GL_{\varphi(k)}$ acts transitively 
 on $Gr(\ell_{\a, k}, V_{\varphi})$. Let 
 $P_{W}$ be the stabilizer of $W$. Then by Serre \cite{S} proposition
 3, we know that the principle bundle 
 $$
 GL_{\varphi(k)}\rightarrow GL_{\varphi(k)}/P_{W}
 $$
 is \'etale-locally trivial.
 Here the base $GL_{\varphi(k)}/P_{W}$ is isomorphic to $Gr(\ell_{\a, k}, V_{\varphi})$.
 It is even Zariski-locally trivial because $P_{W}$ is parabolic, which is special 
 in the sense of Serre \cite{S}, $\S$ 4.
 Now we can write
 \begin{displaymath}
  \xymatrix{
 X_{\a}^{ k}\ar[d] & GL_{\varphi(k)}\times_{P_{W}}\alpha^{-1}(W)\ar[l]_{\hspace{-1.5cm}\delta}\ar[dl]\\
  Gr(\ell_{\a, k}, V_{\varphi})&
  }
\end{displaymath}
 where 
 \[
  \delta([g, T])=g.T.
 \]
We claim that $\delta$ is an isomorphism.
In fact, for any $T\in X_{\a}^{ k}$, 
we choose $g\in GL_{\varphi(k)}$ such that 
\[
 g(\ker(T|_{V_{\varphi, k}}))=W.
\]
This implies $g.T\in \alpha^{-1}(W)$, thus
\[
 \delta([g^{-1}, g.T])=T.
\]
This shows the surjectivity. For injectivity, 
it is enough to show that 
$$
\delta([g,T])=g.T\in \alpha^{-1}(W)
$$
implies $g\in P_{W}$. But this is by definition of $P_{W}$.

The fact that $\alpha$ is locally trivial then can 
be deduced from that of  
$$GL_{\varphi(k)}\times_{P_{W}}\alpha^{-1}(W),$$
while the latter is a consequence of the fact that $GL_{\varphi(k)}$
is locally trivial over $Gr(\ell_{\a, k}, V_{\varphi})$.

Finally, we want to show the surjectivity of the orbit $\alpha|_{O_{\c}}$. 
This is a consequence the fact that $GL_{\varphi(k)}$ acts transitively 
 on $Gr(\ell_{\a, k}, V_{\varphi})$.
\end{proof}

\begin{prop}\label{prop: 3.2.11}
Let $\c\in \tilde{S}(\a)_{k}$. 
The restriction map
\[
 \gamma_{k}: O_{\c}\rightarrow Z^{k, \a}(\c^{(k)})
\]
is surjective.
\end{prop}

\begin{proof}
Let $(T_{0}, W)\in Z^{k, \a}(\c^{(k)})$. 
Consider 
\[
 m=\sharp\{\Delta\in \c: e(\Delta)=k, \deg(\Delta)\geq 2\}\leq \min\{\ell_{\a, k}, \dim(\ker(T_0|_{V_{\varphi, k-1}}))\}.
\]
We choose a splitting $V_{\varphi, k}=W\oplus V_{\varphi, k}/W$ and 
let $T': V_{\varphi, k-1}\rightarrow W$ be a linear morphism of rank $m$.
Finally, we define $T\in \gamma^{-1}_k((T_{0}, W))$ by letting
\[
 T|_{V_{\varphi, k-1}}=T'\oplus T_{0}|_{V_{\varphi, k-1}}, 
\]
\[
T|_{V_{\varphi, k}}=T_{0}|_{V_{\varphi, k}/W}\circ p_{W},
\]
\[
  T|_{V_{\varphi, i}}=T|_{V_{\varphi, i}},
 \text{ for } i\neq k-1, k.
\]
Let 
\[
\{\Delta\in \c: e(\Delta)=k, \deg(\Delta)\geq 2\}=\{\Delta_1, \cdots, \Delta_m\}, \quad b(\Delta_1)\leq \cdots \leq b(\Delta_m).
\]
We denote $W_i=T_0^{[b(\Delta_1), k-1]}(V_{\varphi, b(\Delta_1)})\cap \ker(T_0|_{V_{\varphi, k-1}})$, then
\[
W_1 \subseteq \cdots \subseteq W_r \subseteq \ker(T_0|_{V_{\varphi, k-1}}).
\]
Then we have $T\in O_\c$ if and only if 
\[
\dim(T'(W_i))-\dim(T(W_{i-1}))=\dim(W_i/W_{i-1}), \quad i=1, \cdots, m.
\]
Since such $T'$ always exists, we are done.
\end{proof}

\begin{notation}
 We fix $W\in Gr(\ell_{\a, k}, V_{\varphi})$, and denote 
 $$(X_{\a}^{ k})_{W}, \quad (Z^{k, \a})_{W}$$ the fibers over $W$.
\end{notation}

\begin{prop}\label{prop: 3.3.13}
The fiber $(X_{\a}^{ k})_{W}$ is normal and irreducible as an algebraic variety
over $\C$.
\end{prop}

\begin{proof}
Note that since 
$\tilde{S}(\a)_{k}$ contains a unique minimal element $\c$, 
the variety $X_{\a}^{ k}$ is contained and is open in the 
irreducible variety $\line{O}_{\c}$. Now by 
\cite{Z4} theorem 1, we know that  $X_{\a}^{ k}$ is actually normal.

By proposition \ref{prop: 4.6.10}, we know that $\alpha$ is a fibration between two varieties 
$X_{\a}^{k}$ and $Gr(\ell_{\a, k}, V_{\varphi})$.
The fact that both are normal and irreducible implies that the fiber $(X_{\a}^{ k})_{W}$ is normal 
and irreducible.
\end{proof}

\remk 
Note that by definition, we are allowed to identify 
$(Z^{k, \a})_{W}$ with $Y_{\a^{(k)}}$. This is what 
we do from now on.

\begin{definition}\label{def: 3.3.13}
 We choose a splitting $V_{\varphi, k}=W\oplus V_{\varphi, k}/W$
and denote by $q_{W}: V_{\varphi, k}\rightarrow W$ the projection.
 We define a morphism $\tau_{W}$ 
  $$
  \tau_{W}(T)=((\gamma_{ k})_{W}(T), q_{W}\circ T|_{V_{\varphi, k-1}}).
  $$
  \end{definition}
  
\remk
  Then we have the following commutative diagram
  \begin{displaymath}
  \xymatrix{
  (X_{\a}^{ k})_{W}\ar[r]^{\hspace{-2cm}\tau_{W}}\ar[d]^{(\gamma_{ k})_{W}} 
  &(Z^{k, \a})_{W}\times \Hom(V_{\varphi, k-1}, W)\ar[dl]^{s}\\
  (Z^{k, \a})_{W}
  }
  \end{displaymath}
 where $s$ is the canonical projection.

 \begin{lemma}\label{lem: 4.6.14}
 The morphism $\tau_{W}$ is injective. 
\end{lemma}

\begin{proof}
Note that any $T\in (X_{\a}^{ k})_{W}$ is determined by 
$(\gamma_{ k})_{W}(T)$ and $T|_{V_{\varphi, k-1}}$.
Furthermore, $T|_{V_{\varphi, k-1}}$ is determined by 
$p_{W}\circ T|_{V_{\varphi, k-1}}$ and $q_{W}\circ T|_{V_{\varphi, k-1}}$.
Since $p_{W}\circ T|_{V_{\varphi, k-1}}$ is a component of 
$(\gamma_{ k})_{W}(T)$, 
it is determined by $(\gamma_{ k})_{W}(T)$ and
$q_{W}\circ T|_{V_{\varphi, k-1}}$. This gives us the injectivity.
\end{proof}

\begin{lemma}
Let $\c\in S(\a)_{k}$ such that 
$\c^{(k)}=(\a^{(k)})_{\min}$. Then 
The image of $ O_{\c}\cap (X_{\a}^{k})_{W}$ is open 
in $O_{\c^{(k)}}\times \Hom(V_{\varphi, k-1}, W)$. 
\end{lemma}

\begin{proof}

Let $\c\in S(\a)_{k}$ such that
$\c^{(k)}=(\a^{(k)})_{\min}$.
We shall use the description in proposition 
\ref{prop: 3.2.2}.
We show that the image of 
\[
  O_{\c}\cap (X_{\a}^{k})_{W}
\]
is open in $O_{\c^{(k)}}\times \Hom(V_{\varphi, k-1}, W)$.

Let $T\in (O_{\c})_{W}$. We check case by case:
\begin{description}
 \item [(1)]If $\varphi(k-1)\leq \varphi(k)-\ell_{\a, k}$, 
 the fact $\c^{(k)}=(\a^{(k)})_{\min}$ implies that 
 $T^{(k)}|_{V_{\varphi, k-1}}$ is injective.
 As a consequence
we have $\Im(T|_{V_{\varphi, k-1}})\cap W=0$. Hence for any element $T_{0}\in \Hom(V_{\varphi, k-1}, W)$
, we define $T_{0}\in O_{\c}$, such that 
 \[
  T_{0}|_{V_{\varphi, k-1}}=T_{0}\oplus T^{(k)}|_{V_{\varphi, k-1}},
 \]
which lies in the fiber over $(\gamma_{k})_{W}^{-1}((T^{(k)}, W))$. 
Since by proposition \ref{prop: 3.2.11}, every element in 
$O_{\c^{(k)}}$ comes from some element in $O_{\c}$, 
hence
\[
 \tau_{W}(O_{\c}\cap (X_{\a}^{k})_{W})=O_{\c^{(k)}}\times \Hom(V_{\varphi, k-1}, W),
\]
which is open.
 
\item[(2)]
If $\varphi(k)-\ell_{\a, k}<\varphi(k-1)< \varphi(k)$,
the fact $\c^{ (k)}=(\a^{(k)})_{\min}$ implies that
the morphism 
\[
 T^{ (k)}|_{V_{\varphi, k-1}}
\]
contains a kernel of dimension 
\[
 \varphi(k-1)-\varphi(k)+\ell_{\a, k}.
\]

Our description of $\c$ in proposition \ref{prop: 3.2.2} shows that in this case  
$$
\dim(\Im(T|_{V_{\varphi, k-1}})\cap W)=\varphi(k-1)-\varphi(k)+\ell_{\a, k}.
$$

In this situation, given an element $T_{0}\in \Hom(V_{\varphi, k-1}, W)$ 
 we define $T'\in E_{\varphi}$, such that 
 \[
  T'|_{V_{\varphi, k-1}}=T_{0}\oplus T^{(k)}|_{V_{\varphi, k-1}},
 \]
\[
 T'|_{V_{\varphi, k}}=T^{(k)}|_{V_{\varphi, k}/W}\circ p_{W},
\]
\[
 T'|_{V_{\varphi, i}}=T^{(k)}, \text{ for } i\neq k-1, k.
\]

By construction and proposition 
\ref{prop: 2.2.4}, we know that 
$T'\in O_{\c}$ if and only if $T'|_{V_{\varphi, k-1}}$ is injective, since no segment in $\c$ ends in $k-1$, as described in 
proposition \ref{prop: 3.2.2}. 
And this is equivalent to say
\[
 T_{0}|_{\ker(T^{( k)}|_{V_{\varphi, k-1}})}
\]
is injective. 
This is an open condition, hence $O_{\c}\cap (X_{\a}^{ k})_{W}$ 
is open in $O_{\c^{( k)}}\times \Hom(V_{\varphi, k-1}, W)$.

\item[(3)]If $\varphi(k-1)\geq \varphi(k)$, then by proposition
\ref{prop: 3.2.2}
\[
 \c^{ (k)}=(\a^{( k)})_{\min}
\]
implies 
\[
 \Im(T|_{V_{\varphi, k-1}})\supseteq W.
\]
Recall the notation from proposition \ref{prop: 3.2.2} 
\[
 \a_{0}=\{\Delta_{1}\succeq \cdots \succeq \Delta_{r}\}.
\]
with $r=\varphi(k-1)-\varphi(k)+\ell_{\a, k}$.
Then 
\[
 \c=((\a^{(k)})_{\min}\setminus \a_{0})\cup 
 \{\Delta_{1}^+\succeq \cdots \succeq \Delta_{\ell_{\a, k}}^+\succeq \Delta_{\ell_{\a, k}+1}
 \succeq \cdots \succeq \Delta_{r}\}.
\]

Let $T_{0}\in \Hom(V_{\varphi, k-1},W)$, we define $T'\in E_{\varphi}$
\[
 T'|_{V_{\varphi, k-1}}=T_{0}\oplus T^{(k)}|_{V_{\varphi, k-1}},
 \]
\[
 T'|_{V_{\varphi, k}}=T^{(k)}|_{V_{\varphi, k}/W}\circ p_{W},
\]
\[
 T'|_{V_{\varphi, i}}=T^{(k)}, \text{ for } i\neq k-1, k.
\]
Consider the following flag over $V_{\varphi, k-1}$,
\[
 \ker(T^{(k)}|_{\varphi, k-1})=V_{r}
 \supseteq \cdots \supseteq V_{1}\supseteq V_{0}=0,
\]
where $V_{i}=\Im((T^{(k)})^{\Delta_{i}})\cap \ker(T^{(k)}|_{\varphi, k-1})$, with $i=1, \cdots, r$, where
$T^{[i,j]}$ is the
 composition map: 
 \begin{displaymath}
\xymatrix
{
V_{i}\ar[r]^{\hspace{-0.5cm }T} & V_{i+1}\cdots \ar[r]^{\hspace{0.5cm }T}&V_{j}.\\
}
\end{displaymath}

Now by proposition \ref{prop: 2.2.4}, we know that 
$T'\in O_{\c}$ if and only if 
\[
 \dim(T_{0}(V_{i}))-\dim(T_{0}(V_{i-1}))
 =\dim(V_{i}/V_{i-1}),
\]
for $i=1, \cdots, \ell_{\a, k}$.
In fact, if $V_{i}\neq V_{i-1}$, then 
\[
 \dim(V_{i}/V_{i-1})=\sharp\{j: \Delta_{j}=\Delta_{i}\}.
\]
And by construction, if $i\leq \ell_{\a, k}$, by proposition \ref{prop: 2.2.4}, 
the fact that $\c$ contains $\Delta_{i}^+$ implies that if $T'\in O_{\c}$, 
\[
 \dim(T_{0}(V_{i}))-\dim(T_{0}(V_{i-1}))
 =\dim(V_{i}/V_{i-1}).
\]
The converse holds by the same reason.

Again, this is an open condition, which proves that 
$O_{\c}\cap (X_{\a}^{k})_{W}$ is open in 
$O_{\c^{(k)}}\times \Hom(V_{\varphi, k-1}, W)$. 
\end{description}
 
\end{proof}

\begin{prop}\label{prop: 4.6.14}

The morphism $\tau_{W}$ is an open immersion.

\end{prop}

\begin{proof}
To see that it is open immersion, we shall use Zariski's main theorem.
Since all Schubert varieties are normal, we observe that 
\[
(Z^{k, \a})_{W}\times \Hom(V_{\varphi, k-1}, W)
\]
are normal by theorem 1 of \cite{Z4}. Also, by 
proposition \ref{prop: 3.3.13}, we know that 
$(X_{\a}^{k})_{W}$ is irreducible and normal,
hence $\tau_W$ is an open immersion.
\end{proof}

\begin{prop}\label{prop: 4.6.15}
Let $\c\in \tilde{S}(\a)_{k}$. Then 
$\c\in S(\a)_{k}$ if and only if
\[
 O_{\c}\cap (X_{\a}^{k})_{W}
\]
is open in 
\[
 (O_{\c^{(k)}}\times \Hom(V_{\varphi, k-1}, W)).
\]

\end{prop}

\begin{proof}
We already showed that  
\[
 O_{\c}\cap (X_{\a}^{k})_{W}
\]
is a  sub-variety of 
\[
 O_{\c^{(k)}}\times \Hom(V_{\varphi, k-1}, W).
\]
Moreover, we know that 
\[
  (O_{\c^{(k)}}\times \Hom(V_{\varphi, k-1}, W))\cap (X_{\a}^{k})_{W}
\]
is open in
\[
  O_{\c^{(k)}}\times \Hom(V_{\varphi, k-1}, W)
\]
since $\tau_{W}$ is open.  
Finally, by proposition \ref{prop: 3.2.11}, 
\begin{align*}
 &( O_{\c^{(k)}}\times \Hom(V_{\varphi, k-1}, W))\cap (X_{\a}^{k})_{W}\\
&=\coprod_{\d \in \tilde{S}(\a)_{k}, \d^{(k)}=\c^{(k)}}O_{\d}\cap (X_{\a}^{k})_{W}.
 \end{align*}
The variety $( O_{\c^{(k)}}\times \Hom(V_{\varphi, k-1}, W))\cap (X_{\a}^{k})_{W}$ 
 is irreducible because $( O_{\c^{(k)}}\times \Hom(V_{\varphi, k-1}, W))$ is irreducible, hence the 
 stratification $\coprod_{\d \in \tilde{S}(\a)_{k}, \d^{(k)}=\c^{(k)}}O_{\d}\cap (X_{\a}^{k})_{W}$
 by locally closed sub-varieties can only contain one term which is open, 
 from the point of view of Zariski topology.
Since for any element 
\[
 \d'\in \{\d \in \tilde{S}(\a)_{k}, \d^{(k)}=\c^{(k)}\},
\]
by (d) of lemma \ref{lem: 3.0.7}, we know that 
there exists $\c'\in S(\a)_{k}$ such that $\d'>\c'$. Hence we conclude that
\[
\{\d \in \tilde{S}(\a)_{k}, \d^{(k)}=\c^{(k)}\},
\]
contains a unique minimal element, which lies in $S(\a)_{k}$.
Now our proposition follows. 
\end{proof}

\begin{cor}\label{cor: 4.6.16}
Let $\a$ be a multisegment  and 
$$\c\in S(\a)_{k},$$
then 
\[
 P_{\a, \c}(q)=P_{\a^{(k)}, \c^{(k)}}(q).
 \]
\end{cor}

\begin{proof}
First of all, by proposition \ref{prop: 4.6.10}
and Kunneth formula, we know that 
\[
 \mathcal{H}^{j}(\line{O}_{\c})_{\a}=\mathcal{H}^{j}(\line{O}_{\c}\cap (X_{\a}^{(k)})_{W})_{\a},
\]
the localization being taken at a point in $O_{\a}\cap (X_{\a}^{(k)})_{W}$.
Now by proposition \ref{prop: 4.6.14} and proposition \ref{prop: 4.6.15}
, we may regard $\line{O}_{\c}\cap (X_{\a}^{(k)})_{W}$
as an open subset of $\line{O}_{\c^{(k)}}\times Hom(V_{\varphi, k-1}, W)$, hence
\[
 \mathcal{H}^{j}(\line{O}_{\c}\cap (X_{\a}^{(k)})_{W})_{\a}=
 \mathcal{H}^{j}(\line{O}_{\c^{(k)}}\times Hom(V_{\varphi, k-1}, W))_{\a^{(k)}}
 \]
and Kunneth formula implies that the latter is equal to 
\[
 \mathcal{H}^{j}(\line{O}_{\c^{(k)}})_{\a^{(k)}}.
\]
\end{proof}

\begin{cor}\label{cor: 4.6.17}
Let $\d\in S(\a)$  such that
 \[
  \d^{(k)}= \a^{(k)},
 \]
and 
 \[
  \c\in S(\a)_{k},
 \]
then $\c< \d$, and 
\[
 P_{\a, \c}(q)=P_{\d, \c}(q).
\]
\end{cor}

\begin{proof}
By corollary \ref{cor: 3.1.11}, we know that there exists $\c'\in S(\a)_{k}$
such that 
\[
 \d>\c', ~\c'^{(k)}=\c^{(k)}.
\]
And proposition \ref{prop: 4.6.15} implies $\c'=\c$.
Finally, applying the corollary \ref{cor: 4.6.16} to 
the pairs $\{\a, \c\}$ and $\{\d, \c\}$ yields the 
result.
\end{proof}

\subsection{Some consequeces}
In this section, we draw some conclusions from what we 
have done before, espectially the properties related to $\psi_{k}$.

\begin{prop}\label{cor: 3.2.3}
 The map
 \begin{align*}
 \psi_{k}: S(\a)_{k} &\rightarrow S(\a^{(k)})\\
\c &\mapsto \c^{(k)}
 \end{align*}
is bijective. Moreover, 
\begin{itemize}
 \item  for $\c\in  S(\a)_{k}$
\[
 m(\c,\a)=m(\c^{(k)}, \a^{(k)}).
\]
\item for $\b, \c\in S(\a)_{k}$, we have $\b>\c$ if and only if
$\b^{(k)}> c^{(k)}$.
\end{itemize}
\end{prop}

\begin{proof}
By proposition \ref{prop: 4.6.15},
we know that $\psi_{k}$ is injective. 
Surjectivity is given by proposition \ref{prop: 1.5.10}.

For $\c\in  S(\a)_{k}$,
\[
 m(\c,\a)=m(\c^{(k)}, \a^{(k)})
\]
is by corollary \ref{cor: 4.6.16} by putting $q=1$, 
and applying theorem \ref{teo: 4.1.5}.

Finally, for $\b, \c\in S(\a)_{k}$, if $\b>\c$, then 
$\c\in S(\b^{(k)}, \b)$, and 
by lemma \ref{lem: 3.0.7}, we know that $\b^{(k)}>\c^{(k)}$. 
Reciprocally, if $\b^{(k)}>\c^{(k)}$, by 
proposition \ref{prop: 4.6.15}, we know that $\line{O}_{\b}\subseteq \line{O}_{\c}$, hence 
$\b>\c$. 

\end{proof}

\begin{cor}\label{lem: 2.3.4}
We have
\begin{itemize}
\item 
\addtocounter{theo}{1}
\begin{equation}\label{equ: (2)}
 \pi(\a^{(k)})=
 \sum_{\c\in S(\a)_{k}}m(\c,\a)L_{\c^{(k)}},
\end{equation}
\item let $\b\in S(\a)$ such that $\b$ satisfies the hypothesis $H_{k}(\a)$
and $\b^{(k)}=\a^{(k)}$, then 
\[
 m(\b, \a)=1, ~S(\a)_{k}=S(\b)_{k}.
\]
\end{itemize}
\end{cor}

\begin{proof}

The first part follows from the fact that 
$\psi_{k}$ is bijective and 
$m(\c, \a)=m(\c^{(k)}, \a^{(k)})$.
For the second part of the lemma, we note that 
$L_{\b^{(k)}}=L_{\a^{(k)}}$
appears with multiplicity one in $\pi(\a^{(k)})$, then equation (\ref{equ: (2)})
implies $m(\b, \a)=m(\b^{(k)}, \a^{(k)})=1$.
To see that $ S(\a)_{k}= S( \b)_{k}\subseteq S(\b)$.
Note that we have $S(\b)_{k}\subseteq S(\a)_{k}$ and two bijection 
\[
 \psi_{k}: S(\a)_{k}\rightarrow S(\a^{(k)}),
\]
\[
 \psi_{k}: S(\b)_{k}\rightarrow S(\b^{(k)})=S(\a^{(k)}),
\]
Hence comparing the cardinality gives 
$S(\a)_{k} =S(\b)_k$.
\end{proof}

\section{Reduction to symmetric case}\label{section-reduction-symmetric-case}
\subsection{Minimal Degree Terms}

The goal of this section is to define the 
set $S( \a)_{\d}\subseteq S(\a)$ and describe some of 
its properties.

\begin{definition}
Let  $(k_{1}, \cdots, k_{r})$ be a sequence of integers.
We define 
\[
\a^{(k_{1}, \cdots, k_{r})}=(((\a^{(k_{1})})\cdots)^{ (k_{r})}).
\]
\end{definition}

\begin{notation}
Let $\Delta=[k, \ell]$, we denote 
\[
 \a^{ (\Delta)}=\a^{(k, \cdots, \ell)}.
\]
More generally, for $\d=\{\Delta_{1}\preceq \cdots \preceq \Delta_{r}\}$,
 let
\[
 \a^{ (\d)}=(\cdots ((\a^{ (\Delta_{r})})^{ (\Delta_{r-1})})\cdots)^{ (\Delta_{1})}.
\]
\end{notation}

\begin{definition}
Let $(k_{1}, \cdots, k_{r})$ be a sequence of integers 
, then we define 
\[
S(\a)_{k_{1}, \cdots, k_{r}}=
\{\c\in S(\a):
\c^{( k_{1}, \cdots , k_{i-1})}
\in  S(\a^{( k_{1}, \cdots , k_{i-1})})_{k_{i}}, \text{ for }i=1, \cdots, r\}, 
\]
with the convention 
\[
k_{0}=-\infty, \quad a^{(-\infty)}=a, \quad c^{(-\infty)}=c
\]
and 
\[
 \psi_{k_{1}, \cdots, k_{r}}:
 S(\a)_{k_{1}, \cdots, k_{r}}\rightarrow S(\a^{(k_{1}, \cdots, k_{r})}),
\]
sending $\c$ to $\c^{(k_{1}, \cdots, k_{r})}$.
\end{definition}

\begin{notation} 
Let $\d=\{\Delta_{1}\preceq \cdots \preceq \Delta_{r}\}$ such 
that $\Delta_{i}=[k_{i}, \ell_{i}]$. 
We denote 
\[
 S(\a)_{\d}: =S(\a)_{k_{r},\cdots, \ell_{r}, k_{r-1}, \cdots, k_{1}, \cdots, \ell_{1}}
\]
and 
\[
 \psi_{\d}: =\psi_{k_{r},\cdots, \ell_{r}, k_{r-1}, \cdots, k_{1}, \cdots, \ell_{1}}.
\]
\end{notation}

\begin{prop}\label{prop: 3.2.17}
Let $(k_{1}, \cdots, k_{r})$ be a sequence of integers. Then we have a bijective morphism 
\[
 \psi_{k_{1}, \cdots, k_{r}}: 
 S(\a)_{k_{1}, \cdots, k_{r}}\rightarrow S(\a^{(k_{1}, \cdots, k_{r})}).
 \]
 Moreover,
 \begin{description}
  \item[(1)] For $\c\in S(\a)_{k_{1}, \cdots, k_{r}}$,
  we have 
  \[
   m(\c, \a)=m(\c^{(k_{1}, \cdots, k_{r})}, \a^{(k_{1}, \cdots, k_{r})}).
  \]
\item[(2)] For $\b, \c \in S(\a)_{k_{1}, \cdots, k_{r}}$, 
then $\b>\c$ if and only if $\b^{(k_{1}, \cdots, k_{r})}>\c^{(k_{1}, \cdots, k_{r})}$.
 \item[(3)] We have
 \[
  \pi(\a^{(k_{1}, \cdots, k_{r})})
  =\sum_{\c\in S(\a)_{k_{1}, \cdots, k_{r}}}m(\c, \a)L_{\c^{(k_{1}, \cdots, k_{r})}}.
 \]
\item[(4)] Let $\b\in S(\a)_{k_{1}, \cdots, k_{r}}$ and 
$\b^{(k_{1}, \cdots, k_{r})}=\a^{(k_{1}, \cdots, k_{r})}$, then 
\[
 S(\a)_{k_{1}, \cdots, k_{r}}=S( \b)_{k_{1}, \cdots, k_{r}}.
\]

 \end{description}
 
\end{prop}

\begin{proof}
Injectivity follows from the fact 
\[
 \psi_{k_{1}, \cdots, k_{r}}=\psi_{k_{r}}\circ \psi_{k_{r-1}}\circ \cdots \circ \psi_{k_{1}}.
\]

For surjectivity, let $\d\in S(\a^{(k_{1}, \cdots, k_{r})})$,
we construct $\b$ inductively such that 
$\psi_{k_{1}, \cdots, k_{r}}(\b)=\d$. Let $\a_{r}=\d$, 
assume that we already construct $\a_{i}\in  S(\a^{( k_{1}, \cdots,  k_{i})})_{k_{i+1}}$ satisfying
that 
$$\a_{i}^{( k_{i+1} ,\cdots , k_{j})}\in 
S( \a^{( k_{1}, \cdots , k_{j}}))_{k_{j+1}}$$
for all $i< j\leq r$ and $\a_{i}^{( k_{i+1}, \cdots, k_{r})}=\d$.

Note that by the bijectivity of the morphism 
\[
 \psi_{ k_{i}}: S( \a^{( k_{1}, \cdots , k_{i-1})})_{k_{i}}\rightarrow 
 S(\a^{(k_{1}, \cdots, k_{i})}),
\]
there exists a unique $\a_{i-1}\in S( \a^{( k_{1}, \cdots , k_{i-1})})_{k_{i}}$, such that 
\[
\a_{i-1}^{ (k_{i})}=\a_{i}.
\]
Finally, take $\b=\a_{0}\in S(\a)_{k_{1}, \cdots, k_{r}}$.
We show (1)
by induction on $r$. The case for $r=1$ is by proposition
\ref{cor: 3.2.3}. For general $r$, by induction
\[
 m(\c, \a)=m(\c^{(k_{1}, \cdots, k_{r-1})}, \a^{(k_{1}, \cdots, k_{r-1})}),
\]
and now apply the case $r=1$ to the pair 
$\c^{(k_{1}, \cdots, k_{r-1})}, \a^{(k_{1}, \cdots, k_{r-1})}$
gives 
\[
 m(\c^{(k_{1}, \cdots, k_{r-1})}, \a^{(k_{1}, \cdots, k_{r-1})})=
 m(\c^{(k_{1}, \cdots, k_{r})}, \a^{(k_{1}, \cdots, k_{r})}).
\]
Hence 
\[
 m(\c, \a)=m(\c^{(k_{1}, \cdots, k_{r})}, \a^{(k_{1}, \cdots, k_{r})}).
\]
Also, to show (2), it suffices 
to apply successively the proposition \ref{cor: 3.2.3}.
And (3) follows from the bijectivity of $\psi_{k_{1}, \cdots, k_{r}}$
and (1). 
As for (4), we know by definition, 
\[
 S(\a)_{k_{1}, \cdots, k_{r}}\supseteq 
 S( \b)_{k_{1}, \cdots, k_{r}}.
\]
 We know that any for $\c\in S(\a)_{k_{1}, \cdots, k_{r}}$,
we have $\c^{(k_{1}, \cdots, k_{r})}\leq \b^{(k_{1}, \cdots, k_{r})}$, by (2), 
this implies that $\c\leq \b$. Hence we are done.
\end{proof}

Similarly, we have 
\begin{definition}
Let $(k_{1}, \cdots, k_{r})$ be a sequence of integers, 
then we define 
\[
_{k_{r}, \cdots, k_{1}}S(\a)=
\{\c\in S(\a):
{^{( k_{i}, \cdots , k_{1})}\c}
\in  {_{k_{i+1}}S({^{( k_{i}, \cdots , k_{1})}\a})}, \text{ for }i=1, \cdots, r\}. 
\]
and 
\[
 _{k_{r}, \cdots, k_{1}}\psi:
 {_{k_{r}, \cdots, k_{1}}S}(\a)\rightarrow S(^{(k_{r}, \cdots, k_{1})}\a),
\]
sending $\c$ to $^{(k_{r}, \cdots, k_{1})}\c$.
\end{definition}

\begin{notation} 
Let $\d=\{\Delta_{1},\cdots , \Delta_{r}\}$ such 
that $\Delta_{i}=[k_{i}, \ell_{i}]$ with $k_{1}\leq \cdots \leq k_{r}$
We denote 
\[
 _{\d}S(\a): =_{k_{r},\cdots, \ell_{r}, k_{r-1}, \cdots, k_{1}, \cdots, \ell_{1}}S(\a),
\]
and 
\[
 _{\d}\psi: =_{k_{r},\cdots, \ell_{r}, k_{r-1}, \cdots, k_{1}, \cdots, \ell_{1}}\psi.
\]
\end{notation}

\remk
Let $k_{1}, k_{2}$ be two integers. 
In general, we do not have 
\[
{ _{k_{2}}(S(\a)_{k_{1}})}=(_{k_{2}}S(\a))_{k_{1}}.
\]
For example, let $k_{1}=k_{2}=1$, $\a=\{[1], [2]\}$, then 
\[
 { _{k_{2}}(S(\a)_{k_{1}})}=\{\a\}, ~(_{k_{2}}S(\a))_{k_{1}}=\{[1, 2]\}.
\]

\begin{notation}\label{nota: 4.1.8}
We write for multisegments $\d_{1}, \d_{2}, \a$,
\[
 {_{\d_{2}}S(\a)_{\d_{1}}}: =({_{\d_{2}}S(\a))_{\d_{1}}}, ~
 S(\a)_{\d_{1}, \d_{2}}: =(S(\a)_{\d_{1}})_{\d_{2}}.
\]
and 
\[
 {_{\d_{2}}\psi_{\d_{1}}}: =({_{\d_{2}}\psi)_{\d_{1}}},
 ~\psi_{\d_{1}, \d_{2}}: =(\psi_{\d_{1}})_{\d_{2}}
\] 
And for $\b\in S(\a)$, 
\[
 ^{(\d_{2})}\b^{(\d_{1})}: =(^{\d_{2}}\b)^{(\d_{1})},
 ~\b^{(\d_{1}, \d_{2})}: =(\b^{(\d_{1})})^{(\d_{2})}.
\]

\end{notation}

\subsection{Main Result:symmetrization of multisegments}

Now we return to the main question, i.e., the calculation of 
the coefficient $m(\c, \a)$ for $\c\in S(\a)$. Before we go into the details, 
we describe our strategies: 
\begin{description}
 \item[(i)]Find a symmetric multisegment, denoted by $\a^{\sym}$, 
 such that $L_{\a}$ is the minimal degree term in some partial derivative of 
 $L_{\a^{\sym}}$.
 \item[(ii)]For $\c\in S(\a)$, find $\c^{\sym }\in S(\a^{\sym})$, such that 
 we have $m(\c, \a)=m(\c^{\sym}, \a^{\sym})$.
\end{description}

\begin{prop}\label{cor: 5}
Let $\a$ be any multisegment, then there exists an ordinary multisegment $\b$, and
two  multisegments $\c_{i}, i=1,2$ such that 
\[
 \b\in {_{\c_{2}}S(\b)_{\c_{1}}}, ~ \a={^{(\c_{2})}\b^{(\c_{1})}}
\]
\end{prop}

\begin{proof}
Let $\a=\{\Delta_{1}, \cdots, \Delta_{r}\}$ be such that
\[
 \Delta_{1}\preceq \cdots \preceq \Delta_{r},
\]
and  
\[
 e(\Delta_{1})\leq  \cdots < e(\Delta_{j})=\cdots =e(\Delta_{i})< e(\Delta_{i+1})\leq  \cdots,
\]
such that $\Delta_{j}$ is the smallest multisegment in $\a$ such that $e(\Delta_{j})$ appears 
in $e(\a)$ with multiplicity greater than 1. Let 
$\Delta^{1}=[e(\Delta_{i})+1, \ell]$ be a segment, where $\ell$ is the maximal integer such that for any 
$m$ such that 
$e(\Delta_{i})\leq m\leq \ell-1$,
there is a segment in $\a$ which ends in $m$. Let $\a_{1}$ 
be the multisegment obtained 
by 
replacing $\Delta_{i}$ by $\Delta_{i}^+$,  and all $ \Delta \in \a$ such that $e(\Delta)\in (e(\Delta_{i}), \ell]$ by $\Delta^+$.
 Now we continue the previous
construction with $\a_{1}$ to get $\a_{2}\cdots$, until we get a 
multisegment $\a_{r_{1}}$ such that
$e(\a_{r_{1}})$ contains no segment with multiplicity greater
than 1. Let
\[
 c_{1}=\{\Delta^{1}, \Delta^{2}, \cdots, \Delta^{r_{1}}\}.
\]
Note that by construction, we have
\[
 \Delta^{1}\prec \Delta^{2}\prec \cdots \prec \Delta^{r_{1}}.
\]
And we show that $\a_{r_{1}}\in S(\a_{r_{1}})_{\c_{1}}$.
Note that 
\[
 \a_{i}=\a_{r_1}^{(\Delta^{r_{1}}, \cdots, \Delta^{i+1})},
\]
by induction on $r_1$, we can assume  that $\a_{1}\in S(\a_{r_{1}})_{\Delta^{r_{1}}, \cdots, \Delta^{2}}$
and show that 
$\a\in S(\a_{1})_{\Delta^{1}}$.
We observe that in $\a_{1}$, by construction, 
with the notations above, $\Delta_{j}, \cdots, \Delta_{i-1}$ are the only segments
in $\a_{1}$ that ends in $e(\Delta_{i})$, and $\Delta_{i}^+$ is the only segment
in $\a_{1}$ that ends in $e(\Delta_{i})+1$. Hence we conclude that 
$\a_{1}\in S(\a_{1})_{e(\Delta_{i})+1}$. And for $e(\Delta_{i})+1<m\leq \ell$, we know that 
$\a_{1}^{(e(\Delta_{1})+1, \cdots, m-1)}$ does not contain a segment which ends in $m-1$, hence 
$\a_{1}^{(e(\Delta_{1})+1, \cdots, m-1)}\in S(\a_{1}^{(e(\Delta_{1})+1, \cdots, m-1)})_{m}$.
We are done by putting $m=\ell$.

Now same construction can be applied to show that there exists a multisegment
$\a_{r_{2}}$ such that $b(\a_{r_{2}})$ contains no segment with multiplicity greater
than 1, and 
\[
 c_{2}=\{^{1}\Delta, \cdots, ^{r_{2}}\Delta\},
\]
such that
\[
  \a_{r_{2}}\in {^{\c_{2}}S(\a_{2})}, ~\a_{r_{1}}={^{(\c_{2})}\a_{r_{2}}}
\]
as minimal degree component.

Note that in this way we construct an
ordinary multisegment $\b=\a_{r_{2}}$, 
\[
 \b\in {_{\c_{2}}S(\b)_{\c_{1}}}, ~ \a={^{(\c_{2})}\b^{(\c_{1})}}
\]
\end{proof}

To finish our strategy (i), we are reduced to consider
the case of ordinary multisegments.

\begin{prop}\label{prop: 4.2.2}
Let $\b$ be an ordinary multisegment, then there exists a symmetric multisegment $\b^{\sym}$, and
a multisegment $\c$ such that 
such that 
\[
 \b^{\sym}\in S(\b^{\sym})_{\c},~ \b=\b^{\sym, ~(\c)}.
\]
\end{prop}

\begin{proof}
In general $\b$ is not symmetric, i.e,
we do not have $\min\{e(\Delta): \Delta\in \b\}\geq \max\{b(\Delta): \Delta\in \b\}$.
Let 
\[
 \b=\{\Delta_{1}, \cdots, \Delta_{r}\},  \quad b(\Delta_{1})> \cdots >b( \Delta_{r}).
\]
so that
\[
 b(\Delta_{1})=\max\{b(\Delta_i): i=1, \cdots, r\}.
\]
If $\b$ is not symmetric, let
$\Delta^{1}=[\ell, b(\Delta_{1})-1] $ with $\ell$ maximal satisfying 
that for any $m$ such that $\ell-1\leq m\leq b(\Delta_{1})$, there is a segment in $\b$ starting in $m$.
We construct $\b_{1}$ by replacing every segment $\Delta$ in 
$\b$ ending in $\Delta^{1}$ by $^{+}\Delta$.
Repeat this construction with $\b_{1}$ to get $\b_{2}\cdots $,
until we get $\b^{\sym}=\b_{s}$, which is symmetric.
Let $\c=\{\Delta^{1}, \cdots, \Delta^{s}\}$, then as before, we have 
\[
 \b^{\sym}\in {_{\c}S(\b^{\sym})}, ~ \b= ~^{(\c)}(\b^{\sym}).
\]
\end{proof}

As a corollary, we know that 

\begin{cor}\label{cor: 4.2.3}
For any multisegment $\a$, we can 
find a symmetric multisegment $\a^{\sym}$ and three multisegments $\c_{i}, i=1, 2, 3$, such that 
\[
 \a^{\sym}\in {_{\c_{2}, \c_{3}}S(\a^{\sym})_{\c_{1}}}, ~
 \a={^{(\c_{2}, \c_{3})}\a^{\sym, (\c_{1})}}.
\]
\end{cor}

Now applying proposition \ref{prop: 3.2.17}
\begin{prop}\label{prop: 4.2.4}
The morphism 
\[
{_{\c_{2}, \c_{3}}\psi_{\c_{1}}}: {_{\c_{2}, \c_{3}}S(\a^{\sym})_{\c_{1}}}
 \rightarrow S(\a)
\]
is bijective, and  for $\b\in S(\a)$, there exists a unique $\b^{\sym}\in S(\a^{\sym})$ such that
\[
 m(\b, \a)=m(\b^{\sym}, \a^{\sym}).
\]
\end{prop}

\subsection{Examples }

In this section we shall give some examples to illustrate the idea of 
reduction to symmetric case. 

We first take $\a=\{[1], [2], [2], [3]\}$ to show how to reduce 
a general multisegment to an ordinary multisegment. 
The procedure is showed in the following picture. 

\begin{figure}[!ht]
\centering
\includegraphics{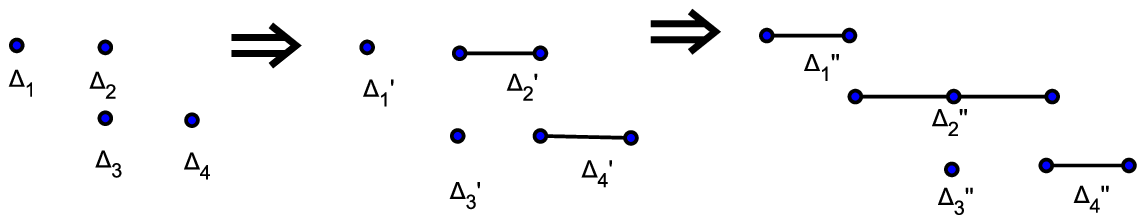}
\caption{\label{fig-multisegment4} }
\end{figure} 
Here we have $\a_{2}=\{[0, 1], [1,3], [2], [3, 4]\}$, 
such that 

\[
 \a_{2}\in {_{[0,1]}S(\a_{2})_{[3,4]}}, ~\a={^{([0,1])}\a_{2}^{([3,4])}}
\]

Next, we reduce the ordinary multisegment $\a_{2}$ to 
a multisegment $\a^{\sym}$, as is showed in the 
following picture.

\begin{figure}[!ht]
\centering
\includegraphics{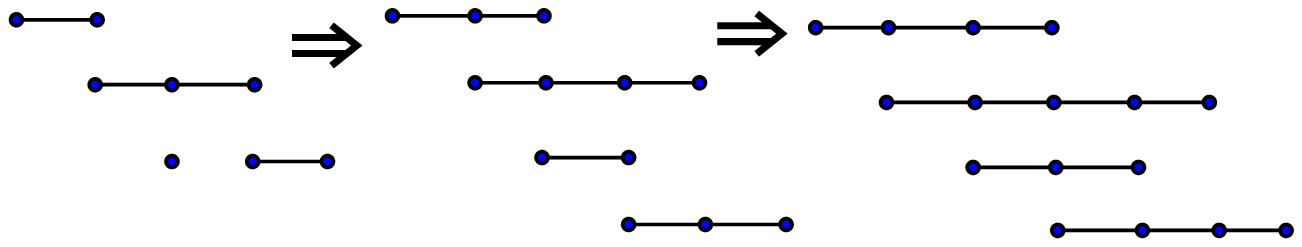}
\caption{\label{fig-multisegment5} }
\end{figure} 

Here,we have 
\[
 \a^{\sym}=\{[0, 3], [1, 5], [2, 4], [3, 6]\}=\Phi(w)
\]
where $w=\sigma_{2}\in S_{4}$.

Now we take $\b=\{[1, 2], [2, 3]\}$, we want to find 
$\b^{\sym}\in S(\a^{\sym})$ such that $m(\b, \a)=m(\b^{\sym}, \a^{\sym})$. 
Actually, following the procedure in Figure 2 above, we have 

\begin{figure}[!ht]
\centering
\includegraphics{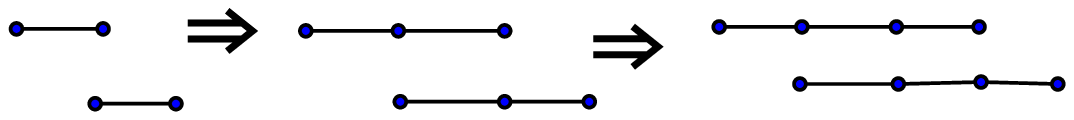}
\caption{\label{fig-multisegment6} }
\end{figure} 

Here we get $\b_{2}=\{[0, 3], [1], [2, 4]\}$. Again, follow the procedure
in Figure 3 above gives

\begin{figure}[!ht]
\centering
\includegraphics{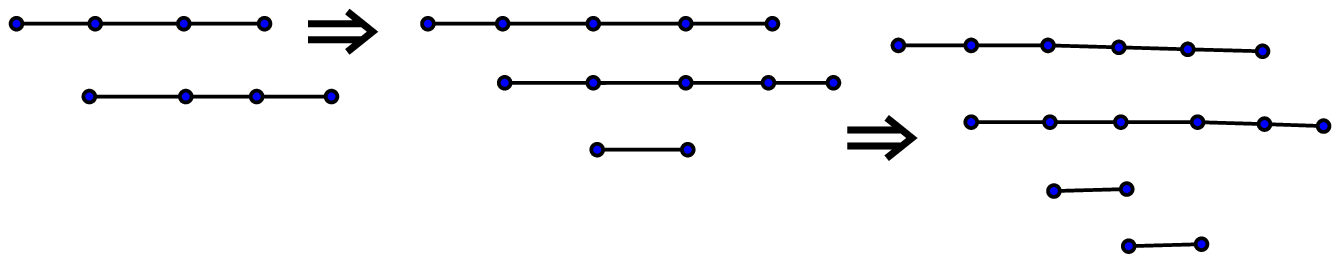}
\caption{\label{fig-multisegment7} }
\end{figure} 

Hence we get  
\[
 \b^{\sym}=\{[0, 5], [1, 3], [2, 6], [3, 4]\}=\Phi(v)
\]
with $v=(13)(24)\in S_{4}$. From \cite{Z2} section 11.3, we know
that $m(\b, \a)=2$, hence we get $m(\b^{\sym}, \a^{\sym})=2$.

\remk
We showed in section 2 that 
\[
 m(\b^{\sym}, \a^{\sym})=P_{v, w}(1),
\]
where $P_{v,w}(q)$ is the Kazhdan Lusztig polynomial associated to $v, w$.
One knows that $P_{v,w}(q)=1+q$, hence $P_{v,w}(1)=2$.

As we have seen, to each multisegment, we have (at least) two different 
ways to attach a Kazhdan Lusztig polynomial:

(1)To use the Zelevinsky construction as described in 
section 4.2. \\
(2)To first construct an associated symmetric multisegment, 
and then attach the corresponding Kazhdan Lusztig polynomial. 

\vspace{0.5cm}
\remk 
In general, for $\a>\b$, (1) gives a polynomial $P_{\a, \b}^{Z}$ which is a
Kazhdan Lusztig polynomial for the symmetric group $S_{\deg(\a)}$. And (2)
gives a polynomial $P_{\a, \b}^{S}$, which is a KL polynomial for a symmetric 
group $S_{n}$ with $n\leq \deg(\a)$. 
It may happen that $n=\deg(\a)$.
By corollary \ref{cor: 4.6.16}, we always have $P_{\a, \b}^{Z}=P_{\a, \b}^{S}$.

\begin{example}
Consider $\a=\{1, 2, 2, 3\}, \b=\{[1,2], [3,4]\}$, then by 
\cite{Z3} section 3.4, we know that $P_{\a, \b}^{Z}=1+q$.
And the symmetrization of $\a$ and $\b$ are given by
\[
 \a^{\sym}=\Psi((2,3)), \quad \b^{\sym}=\Psi((1,3)(2,4)).
\]
Hence $P_{\a, \b}^{S}=P_{(2,3), (1,3)(2,4)}=1+q$,
which is 
the Kazhdan Lusztig polynomial for the pair $((2,3), (1,3)(2,4))$ in $S_{4}$
\end{example}

\section{Application: proof of the Zelevinsky's conjecture}
 \begin{definition}
 The relation type between 2 segments $\{\Delta, \Delta'\}$ is one of the following
 \begin{itemize}
  \item $\Delta$ cover $\Delta'$ if $\Delta\supseteq \Delta'$;
 \item linked  but not juxtaposed if $\Delta$ does not cover $\Delta'$ and $\Delta\cup \Delta'$ is a segment but $\Delta\cap \Delta'\neq \emptyset$;
 \item juxtaposed if $\Delta\cup \Delta'$ is a segment but $\Delta\cap \Delta'=\emptyset$;
 \item unrelated if $\Delta\cap \Delta'=\emptyset$ and $\Delta, \Delta'$ are not linked.
 \end{itemize}

 \end{definition}

\begin{definition}
Two multisegments 
$$\a=\{\Delta_1,  \cdots,  \Delta_r\} \qquad \hbox{and} \qquad 
\a'=\{\Delta_1', \cdots, \Delta'_{r'}\}$$ 
have the same relation type if 
\begin{itemize}
\item $r=r'$;
\item there exists a bijection 
\[
\xi: \a\rightarrow \a'
\]
 of multisets which preserves the partial order $\preceq$ and relation type of segments and induces bijection of multisets 
 \[
 e(\xi): e(\a)\rightarrow e(\a'), \quad b(\xi): b(\a)\rightarrow b(\a'). 
 \]
\end{itemize} 
satisfying 
\[
e(\xi)(e(\Delta))=e(\xi(\Delta)), \quad b(\xi)(b(\Delta))=b(\xi(\Delta)).
\]
\end{definition}

\begin{lemma}
Let $\a$ and $\a'$ be of the same relation type induced by $\xi$. 
Let $\{\Delta_1\preceq \Delta_2\}$ be linked in $\a$.  Denote by $\a_1$($\a_1'$, resp.) the multisegment obtained by applying the elementary
operation to $\{\Delta_1, \Delta_2\}$( $\{\xi(\Delta_1), \xi(\Delta_2)\}$, resp.). 
Then $\a_1$ and $\a_1'$
also have the same relation type.
\end{lemma}
\begin{proof}
We define a bijection
\[
\xi_1: \a_1\rightarrow \a_1'
\]
by
\[
\xi_1(\Delta_1\cup \Delta_2)=\xi(\Delta_1)\cup \xi(\Delta_2), \quad \xi_1(\Delta_1\cap \Delta_2)=\xi(\Delta_1)\cap \xi(\Delta_2)
\]
and 
\[
\xi_1(\Delta)=\xi(\Delta), \quad \text{ for all } \Delta\in \a\setminus \{\Delta_1, \Delta_2\}.
\]
It induces a bijection between the end multisets $e(\a_1)$ and $e(\a_1')$ as well as the beginning multisets 
$b(\a_1)$  and $b(\a_1')$.
Also the morphism $\xi$ preserves the partial order follows from the fact that for 
$x, y\in e(\a)$ such that $x\leq y$,  then $e(\xi_1)(x)=e(\xi)(x)\leq e(\xi_1)(y)=e(\xi)(y)$( The same fact holds for  $b(\xi_1)$).
Finally, it remains to show that $\xi_1$ respects the relation type. 
Let $\Delta\preceq \Delta'$ be two segments in $\a_1$, if non of them is contained in $\{\Delta_1\cup \Delta_2, \Delta_1\cap \Delta_2\}$, 
then $\xi_1(\Delta)=\xi(\Delta)$ and $\xi_1(\Delta')=\xi(\Delta')$  and they are in the same relation type as $\{\Delta, \Delta'\}$ by assumption.
For simplicity, we only discuss the case where $\Delta=\Delta_1\cup \Delta_2$ 
but $\Delta'$ is not contained in $\{\Delta_1\cup \Delta_2, \Delta_1\cap \Delta_2\}$, other cases are similar. 
\begin{itemize}
\item If $\Delta'$ cover $\Delta$, then $\Delta$ cover $\Delta_1$ and $\Delta_2$,  hence
$\xi_1(\Delta)=\xi(\Delta)$ cover $\xi(\Delta_1)$ and $\xi(\Delta_2)$, which implies $\xi_1(\Delta')$ covers $\xi_1(\Delta)$.
\item If $\Delta'$ is linked to $\Delta$ but not juxtaposed, then either $\Delta'$ covers $\Delta_2$ and linked to $\Delta_1$, 
or $\Delta'$ is linked to $\Delta_2$ but not juxtaposed. In both cases we have $\xi(\Delta')$ is linked to $\xi(\Delta_1)\cup \xi(\Delta_2)$ and 
not juxtaposed. 
\item If $\Delta'$ is juxtaposed to $\Delta$, then $\Delta'$ is juxtaposed to $\Delta_2$ since $\Delta_2\succeq \Delta_1$. 
Therefore $\xi(\Delta')$ is juxtaposed to $\xi(\Delta_2)$ which implies $\xi_1(\Delta')$ is juxtaposed to the segment $\xi_1(\Delta)$.
\item If $\Delta'$ is unrelated to $\Delta_1\cup \Delta_2$, then it is unrelated to both $\Delta_1$ and $\Delta_2$ with $\Delta_2\preceq \Delta'$, 
this implies that $\xi(\Delta')$ is unrelated to $\xi(\Delta_1)\cup \xi(\Delta_2)$.  
\end{itemize}

\end{proof}

\remk As every element $\b\in S(\a)$ is obtained from $\a$ by a sequence of elementary operations,
we can define an application of poset
$$
\Xi: S(\a)\longrightarrow S(\a').
$$

\begin{lemma}
 The application $\Xi$ is well defined and bijective.
 \end{lemma}
 \begin{proof}
 We give a new definition of $\Xi$ in the following way. For $\b\in S(\a)$, we define
 \[
 \Xi(\b)=\{[b(\xi)(b(\Delta)), e(\xi)(e(\Delta))]: \Delta\in \b\}
 \]
 such a definition is independent of the choice of elementary operations.
 It remains to see that it  coincides with the one using elementary operation. 
 In fact, let $\a_1$ be a multisegment obtained by applying the elementary operation to 
 the pair of segments $\{\Delta_1\preceq \Delta_2\}$, then by our original definition of $\Xi$, it
sends $\a_1$ to $\a_1'$ in the previous lemma. Now 
 by the new definition, we have  $\Xi(\a_1)$ given by
 \[
\{\xi(\Delta):   \Delta\in \a\setminus \{\Delta_1, \Delta_2\}\}\cup \{[b(\xi)(b(\Delta_1)), b( \xi)(b(\Delta_2))], [b(\xi)(b(\Delta_2)), b( \xi)(b(\Delta_1))]\}.
 \]
 By our definition of $\xi$, we get 
 \[
 [b(\xi)(b(\Delta_1)), b( \xi)(b(\Delta_2))]=\xi(\Delta_1)\cup \xi(\Delta_2),
 \]
 and 
 \[
 [b(\xi)(b(\Delta_2)), b( \xi)(b(\Delta_1))]=\xi(\Delta_1)\cap \xi(\Delta_2).
 \]
 
 Hence we conclude that $\Xi$ is well defined.  
 Note that by our definition, $\xi$ is invertible, which gives $\xi^{-1}$, and in the same way we can construct $\Xi^{-1}$.
 Now we have 
 \[
 \Xi\Xi^{-1}=\Id, \quad \Xi^{-1}\Xi=\Id
 \]
 by our definition above using $b(\xi)$ and $e(\xi)$.  This shows that $\Xi$ is bijective.

 \end{proof}

%
%

\begin{teo}\label{teo: 4.4.5}
For $\a$ and $\a'$ having the same relation type,  then for $\b\in S(\a)$ with $\b'=\Xi(\b)$, we have 
$$
m(\b, \a)=m(\b', \a').
$$
\end{teo}

\begin{proof}
First of all, we consider the case where $\a$ and $\a'$ are symmetric multisegments.
Let $\a=\Phi(w)$ by fixing a map
\[
\Phi: S_n\rightarrow S(\a_{\Id}).
\]
Now since $\a$ and $\a'$ have the same relation type, we know that $\a'=\Phi'(w)$
for some fixe map 
\[
\Phi': S_n\rightarrow S(\a'_{\Id}).
\]
Finally,  let $\a=\{\Delta_1, \cdots, \Delta_n\}$ and $\a'=\{\Delta_1', \cdots, \Delta_n'\}$ such that 
\[
b(\Delta_1)<\cdots< b(\Delta_n),  \quad \Delta_i'=\xi(\Delta_i).
\]
Without loss of generality, we assume that $b(\Delta_1)=b(\Delta_1')$. 
We can assume that $b(\Delta_i)=b(\Delta_{i-1})+1$. In fact, if 
$b(\Delta_i)>b(\Delta_{i-1})+1$, then by replacing $\Delta_i$ by $^{+}\Delta_i$ , we get a new symmetric multisegment $\a_1$
which has the same relation type as $\a$. Moreover, let $\b\in S(\a)$ and $\b_1$ be the corresponding multisegment in $S(\a_1)$, then 
\[
m(\b, \a)=m(\b_1, \a_1)
\]
by proposition \ref{cor: 3.2.3}. 
We note that the equality
\[
m(\b_1, \a_1)=m(\b', \a')
\]
implies that
\[
m(\b', \a')=m(\b, \a).
\]
Therefore it suffices to prove the theorem for $\a_1$ and $\a'$.
From now on, let $b(\Delta_i)=b(\Delta_{i-1})+1$ and $b(\Delta_i)=b(\Delta_i')$. 
The same argument shows that we can furthermore assume that 
\[
e(\Delta_{w^{-1}(i)})=e(\Delta_{w^{-1}(i-1)})+1, \quad e(\Delta'_{w^{-1}(i)})=e(\Delta'_{w^{-1}(i-1)})+1.
\]
Now if $e(\Delta_{w^{-1}(1)})<e(\Delta'_{w^{-1}(1)})$, 
then consider the truncation functor  $\a'\mapsto \a'^{(e(\Delta_{w^{-1}(1)})+1, \cdots, e(\Delta_{w^{-1}(1)}))}$, 
the latter is a symmetric multisegment having the same relation type as $\a'$, and 
\[
m(\b', \a')=m(\b'^{(e(\Delta_{w^{-1}(1)})+1, \cdots, e(\Delta_{w^{-1}(1)}))}, \a'^{(e(\Delta_{w^{-1}(1)})+1, \cdots, e(\Delta_{w^{-1}(1)}))})
\] 
by proposition \ref{prop: 3.2.17}. Repeat the same procedure, 
in finite step, we find $\c$, such that 
\[
\a=\a'^{(\c)}
\]
and 
\[
m(\b, \a)=m(\b', \a'). 
\]
by  proposition \ref{prop: 3.2.17}.

\remk an interesting application of this computation is given in the corollary \ref{coro-sym1}.

For general case, note that in section 4.4, we construct a symmetric multisegment $\a^{\sym}$  and three multisegments $\c_i, i=1, 2, 3$
such that 
\[
 \a^{\sym}\in {_{\c_{2}, \c_{3}}S(\a^{\sym})_{\c_{1}}}, ~
 \a={^{(\c_{2}, \c_{3})}\a^{\sym, (\c_{1})}}.
\]
(cf. Corollary \ref{cor: 4.2.3}).  
The same for $\a'$, we have 
\[
 \a'^{\sym}\in {_{\c'_{2}, \c'_{3}}S(\a'^{\sym})_{\c'_{1}}}, ~
 \a'={^{(\c'_{2}, \c'_{3})}\a'^{\sym, (\c'_{1})}}.
\]

\end{proof}

\begin{lemma}
The two multisegment $\a^{\sym}$ and $\a'^{\sym}$ have the same relation type. And let 
$\Xi^{\sym}: S(\a^{\sym})\rightarrow \a'^{\sym}$ be the bijection constructed above, then we have 
the following commutative diagram
\begin{displaymath}
\xymatrix{
{_{\c_{2}, \c_{3}}S(\a^{\sym})_{\c_{1}}}\ar[r]^{\Xi^{\sym}}\ar[d]^{{_{\c_{2}, \c_{3}}\psi_{\c_{1}}}}& {_{\c'_{2}, \c'_{3}}S(\a'^{\sym})_{\c'_{1}}}\ar[d]^{_{\c'_{2}, \c'_{3}}\psi_{\c'_{1}}}\\
S(\a)\ar[r]^{\Xi}& S(\a').
}
\end{displaymath}
\end{lemma}

Admitting the lemma, we have 
\[
m(\b, \a)=m(\b^{\sym}, \a^{\sym}), \quad m(\b', \a')=m(\b'^{\sym}, \a'^{\sym})
\]
by proposition \ref{prop: 4.2.4}. Now by what we have proved before and the above lemma, we have 
\[
m(\b^{\sym}, \a^{\sym})=m(\b'^{\sym}, \a'^{\sym}),
\]
which implies $m(\b, \a)=m(\b', \a')$.

\begin{proof}
Note that by construction we know that the number of segments in $\a^{\sym}$ is the same as that of $\a$. 
Let $\a^{\sym}=\{\Delta_1\preceq \cdots \preceq  \Delta_r\}$, then 
$\a=\{{^{(\c_{2}, \c_{3})}\Delta_1^{ (\c_{1})}}\preceq \cdots \preceq {^{(\c_{2}, \c_{3})}\Delta_r^{(\c_{1})}}\}$.
Also let $\a'^{\sym}=\{\Delta'_1\preceq \cdots \preceq  \Delta'_r\}$.
We define 
\begin{align*}
\xi^{\sym}:& \a^{\sym}\rightarrow \a'^{\sym}\\
  & \Delta_i \mapsto\Delta'_i.
\end{align*}
This automatically induces bijections
\[
e(\xi^{\sym}): e(\a^{\sym})\rightarrow e(\a'^{\sym}), \quad b(\xi^{\sym}): b(\a^{\sym})\rightarrow b(\a'^{\sym}),
\]
since all of them are sets.  Note that we definitely have 
\[
\xi({^{(\c_{2}, \c_{3})}\Delta_i^{ (\c_{1})}})={^{(\c'_{2}, \c'_{3})}\Delta_i'^{ (\c'_{1})}}
\]

It remains to show that $\xi^{\sym}$ preserve the relation type.  
Let $i\leq j$.  Then $\Delta_i$ and $\Delta_j$ are linked if and only if one of the following happens 
\begin{itemize}
\item ${^{(\c_{2}, \c_{3})}\Delta_i^{ (\c_{1})}}$ and ${^{(\c_{2}, \c_{3})}\Delta_j^{ (\c_{1})}} $ are linked, juxtaposed or not;
\item ${^{(\c_{2}, \c_{3})}\Delta_i^{ (\c_{1})}} $ and ${^{(\c_{2}, \c_{3})}\Delta_j^{ (\c_{1})}}$ are unrelated.
\end{itemize}

And $\Delta_j$ covers $\Delta_i$ if and only if ${^{(\c_{2}, \c_{3})}\Delta_j^{ (\c_{1})}}$ covers ${^{(\c_{2}, \c_{3})}\Delta_i^{ (\c_{1})}}$. 
Since $\xi$ preserves relation types, this shows that $\xi^{\sym}$ also preserves relation types.  Hence we conclude that 
$\a^{\sym}$ and $\a'^{\sym}$ have same relation type.  To see that the map $\Xi^{\sym}$ sends ${_{\c_{2}, \c_{3}}S(\a^{\sym})_{\c_{1}}}$ to 
${_{\c'_{2}, \c'_{3}}S(\a'^{\sym})_{\c'_{1}}}$, consider $\b\in S(\a)$ and its related element $\b^{\sym}\in {_{\c_{2}, \c_{3}}S(\a^{\sym})_{\c_{1}}}$. 

\medskip

- First of all, we assume that $l(\b)=1$, i.e. $\b$ can be obtained from $\a$ by applying the elementary operation to the pair 
$\{{^{(\c_{2}, \c_{3})}\Delta_i^{ (\c_{1})}}, {^{(\c_{2}, \c_{3})}\Delta_j^{ (\c_{1})}}\}(i<j)$.  Let 
$\tilde{\b}$ be the element in $S(\a^{\sym})$ obtained by applying the elementary operation to the pair of segments $\{\Delta_i, \Delta_j\}$ in $\a^{\sym}$. 
Then we have 
\[
\b={^{(\c_{2}, \c_{3})}\tilde{\b}^{(\c_{1})}}.
\]
Let $\tilde{\b}'=\Xi^{\sym}(\tilde{\b})$.  By construction, we have
\[
\b'=\Xi(\b)={^{(\c'_{2}, \c'_{3})}\tilde{\b}'^{(\c'_{1})}}.
\]
Now consider 
\[
\tilde{\b}_0=\tilde{\b}>\cdots >\tilde{\b}_n=\b^{\sym}
\]
be a maximal chain of multisegments and let $\tilde{\b}'_i=\Xi^{\sym}(\tilde{\b}'_i)$, then  
\[
\tilde{\b}'_0>\cdots >\tilde{\b}'_n.
\]
Let 
\[
\tilde{\b}_{i}=\{\Delta_{i, 1}\preceq \cdots \preceq \Delta_{i, r_i}\},\quad \tilde{\b}'_{i}=\{\Delta_{i, 1}'\preceq \cdots \preceq \Delta_{i, r_i}'\}.
\]
We prove by induction that 
\[
\b'={^{(\c'_{2}, \c'_{3})}\tilde{\b}_i'^{(\c'_{1})}}.
\]
We already showed the case where $i=0$. Assume that we have 
\[
\b'={^{(\c'_{2}, \c'_{3})}\tilde{\b}_j'^{(\c'_{1})}}
\]
for $j<i$.  Suppose that $\tilde{\b}_i$ is obtained from $\tilde{\b}_{i-1}$ by applying the elementary operation to 
the pair of segments $\{\Delta_{i-1, \alpha_{i-1}}\preceq \Delta_{i-1, \beta_{i-1}}\}$.  
We deduce from the fact $\tilde{\b}_i\geq \b^{\sym}$ that we are in one of the following situatios
\begin{itemize}
\item  ${^{(\c_{2}, \c_{3})}\Delta_{i-1, \alpha_{i-1}}^{(\c_{1})}}=\emptyset$ or ${^{(\c_{2}, \c_{3})}\Delta_{i-1, \beta_{i-1}}^{(\c_{1})}}=\emptyset$;
\item  $b({^{(\c_{2}, \c_{3})}\Delta_{i-1, \beta_{i-1}}^{(\c_{1})}})=b({^{(\c_{2}, \c_{3})}\Delta_{i-1, \alpha_{i-1}}^{(\c_{1})}})$;
\item $e({^{(\c_{2}, \c_{3})}\Delta_{i-1, \beta_{i-1}}^{(\c_{1})}})=e({^{(\c_{2}, \c_{3})}\Delta_{i-1, \alpha_{i-1}}^{(\c_{1})}})$.
\end{itemize}
According the our assumption that $\tilde{\b}'_i=\Xi^{\sym}(\tilde{\b}'_i)$, we have 
\[
\xi({^{(\c_{2}, \c_{3})}\Delta_{i-1, j}^{(\c_{1})}})={^{(\c_{2}, \c_{3})}\Delta_{i-1, j}'^{(\c_{1})}},
\]
therefore the pair $\{{^{(\c_{2}, \c_{3})}\Delta_{i-1, \alpha_{i-1}}'^{(\c_{1})}}, {^{(\c_{2}, \c_{3})}\Delta_{i-1, \beta_{i-1}}'^{(\c_{1})}}\}$
also satisfies one of the listed properties above.  And this shows that $\tilde{\b}'_i$ is sent to $\b'$ by $_{\c'_{2}, \c'_{3}}\psi_{\c'_{1}}$. 
Therefore by proposition \ref{prop: 4.6.15}, we know that 
\[
\b_n'\geq \b'^{\sym}.
\]
Conversely, we have 
\[
\Xi^{\sym-1}(\b'^{\sym})\geq \b^{\sym}.
\]
Combine the two inequalities to get
\[
\Xi^{\sym}(\b^{\sym})= \b'^{\sym}.
\]

- The general case where $\ell(\b)>1$, we can choose a maximal chain of multisegments 
 \[
 \a=\a_0>\cdots >\a_{\ell(\b)}=\b.
 \]
 Let $\a_i'=\Xi(\a_i)$, by assumption, we can assume that for $i<\ell(\b)$, we have
 \[
 \Xi^{\sym}(\a_i^{\sym})=\a_i'^{\sym}.
 \]
 By considering the set $S(\a_{\ell(\b)-1})$, we are reduce to the case where $\ell(\b)=1$. Hence we are done.
 \end{proof}
 
\begin{cor} \label{coro-sym1}
Let $\a_{\Id}$ be a symmetric multisegment associated to the identity in $S_n$ and 
\[
\Phi: S_n\rightarrow S(\a_{\Id}).
\] 
Then 
\[
m(\Phi(v), \Phi(w))=P_{w, v}(1).
\]
\end{cor} 
\begin{proof}
The special case where
\[
\a_{\Id}=\sum_{i=1}^{n}[i, i+n-1]
\]
is already treated in corollary \ref{cor: 2.5.9}.
The general case can be deduced from the theorem above.
\end{proof} 
 
\bibliographystyle{plain}
\bibliography{biblio}

\def\cprime{$'$}
\begin{thebibliography}{10}

\bibitem{Z1}
I.~N. Bernstein and A.~V. Zelevinsky.
\newblock Induced representations of reductive {$p$}-adic groups. {I}.
\newblock {\em Ann. Sci. \'Ecole Norm. Sup. (4)}, 10(4):441--472, 1977.

\bibitem{BF}
Anders Bj{\"o}rner and Francesco Brenti.
\newblock {\em Combinatorics of {C}oxeter groups}, volume 231 of {\em Graduate
  Texts in Mathematics}.
\newblock Springer, New York, 2005.

\bibitem{CG}
Neil Chriss and Victor Ginzburg.
\newblock {\em Representation theory and complex geometry}.
\newblock Modern Birkh\"auser Classics. Birkh\"auser Boston, Inc., Boston, MA,
  2010.
\newblock Reprint of the 1997 edition.

\bibitem{LM17}
A.~M{\'i}nguez E.~Lapid.
\newblock Geometric conditions for square-irreducibility of certain
  representations of the general linear group over a non-archimedean local
  field.
\newblock {\em preprint}.

\bibitem{GLS}
Christof Gei{\ss}, Bernard Leclerc, and Jan Schr{\"o}er.
\newblock Kac--moody groups and cluster algebras.
\newblock {\em Advances in Mathematics}, 228(1):329--433, 2011.

\bibitem{KKMO}
Seok-Jin Kang, Masaki Kashiwara, Myungho Kim, and Se-Jin Oh.
\newblock Monoidal categorification of cluster algebras ii.
\newblock {\em arXiv preprint arXiv:1502.06714}, 2015.

\bibitem{KL}
David Kazhdan and George Lusztig.
\newblock Schubert varieties and {P}oincar\'e duality.
\newblock In {\em Geometry of the {L}aplace operator ({P}roc. {S}ympos. {P}ure
  {M}ath., {U}niv. {H}awaii, {H}onolulu, {H}awaii, 1979)}, Proc. Sympos. Pure
  Math., XXXVI, pages 185--203. Amer. Math. Soc., Providence, R.I., 1980.

\bibitem{MV}
Alberto M{\'{\i}}nguez and Vincent S{\'e}cherre.
\newblock L'involution de zelevinsky modulo $\ell$.
\newblock {\em Preprint}, 2015.

\bibitem{S}
Jean-Pierre Serre.
\newblock Espaces fibr\'es alg\'ebriques (d'apr\`es {A}ndr\'e {W}eil).
\newblock In {\em S\'eminaire {B}ourbaki, {V}ol.\ 2}, pages Exp.\ No.\ 82,
  305--311. Soc. Math. France, Paris, 1995.

\bibitem{Z2}
A.~V. Zelevinsky.
\newblock Induced representations of reductive {$p$}-adic groups. {II}. {O}n
  irreducible representations of {${\rm GL}(n)$}.
\newblock {\em Ann. Sci. \'Ecole Norm. Sup. (4)}, 13(2):165--210, 1980.

\bibitem{Z3}
A.~V. Zelevinsky.
\newblock A $p$-adic analog of the {K}azhdan-{L}usztig conjecture.
\newblock {\em Funct.Anal.Appl.}, 15:83--92, 1981.

\bibitem{Z4}
A.~V. Zelevinsky.
\newblock Two remarks on graded nilpotent classes.
\newblock {\em Uspekhi Mat. Nauk}, 40(1(241)):199--200, 1985.

\end{thebibliography}

\end{document}